\documentclass[preprint,12pt]{elsarticle}

\usepackage{graphicx}
\usepackage{amssymb}
\usepackage[utf8]{inputenc}
\usepackage{mathtools,enumitem,amsmath}
\usepackage{graphicx,xcolor,graphbox}
\usepackage[margin=.5in]{caption}
\usepackage{subcaption}
\usepackage[margin=1in]{geometry}
\usepackage{amsthm}
\usepackage[colorlinks,citecolor=blue]{hyperref}
\usepackage{bm}
\usepackage{here}
\usepackage{booktabs}

\newcommand{\N}{\ensuremath{\mathbb{N}}}
\newcommand{\M}{\ensuremath{\mathbb{M}}}
\newcommand{\Bh}{\ensuremath{\mathcal{B}^{(h)}}}
\newcommand{\G}{\ensuremath{\mathcal{G}}}
\newcommand{\R}{\ensuremath{\mathbb{R}}}
\newcommand{\PN}{\ensuremath{P_{\N^+}}}
\newcommand{\xG}{{\bf \xi}_{\Gamma}}
\newcommand{\xGs}{{\bf \xi}_{\Gamma^*}}
\newcommand{\xg}{{\bf \xi}_{\gamma}}
\newcommand{\uh}{\ensuremath{u^{(h)}}}
\newcommand{\uhn}{\ensuremath{\uh_n}}
\newcommand{\uhnh}{\ensuremath{\uh_{\frac{n}{2}}}}
\newcommand{\selfnorm}{\ensuremath{\|\uhn - \uhnh\|_{\infty}}}
\newcommand{\Ex}{\ensuremath{{\bf Ex}}}
\renewcommand{\L}{\ensuremath{\mathcal{L}}}
\newcommand{\Lh}{\ensuremath{\mathcal{L}^{(h)}}}
\newcommand{\Gh}{\ensuremath{\mathcal{G}^{(h)}}}
\newcommand \defeq{\stackrel{\mathclap{\normalfont\mbox{\tiny def}}}{=}}
\newcommand{\dO}{\partial \Omega}
\newcommand \partialdiff[2]{\ensuremath{\frac{\partial #1}{\partial #2}}}
\journal{Journal of Computational and Applied Mathematics}

\begin{document}

\begin{frontmatter}

%% Title, authors and addresses

\title{Non-Iterative Domain Decomposition for the Helmholtz Equation Using the Method of Difference Potentials\tnoteref{title1}}

\tnotetext[title1]{Work supported by the US Army Research Office (ARO) under grant
W911NF-16-1-0115 and  the US--Israel Binational Science Foundation (BSF) under grant 2014048.}
%% use the tnoteref command within \title for footnotes;
%% use the tnotetext command for the associated footnote;
%% use the fnref command within \author or \address for footnotes;
%% use the fntext command for the associated footnote;
%% use the corref command within \author for corresponding author footnotes;
%% use the cortext command for the associated footnote;
%% use the ead command for the email address,
%% and the form \ead[url] for the home page:
%%
%% \title{Title\tnoteref{label1}}
%% \tnotetext[label1]{}
%% \author{Name\corref{cor1}\fnref{label2}}
%% \ead{email address}
%% \ead[url]{home page}
%% \fntext[label2]{}
%% \cortext[cor1]{}
%% \address{Address\fnref{label3}}
%% \fntext[label3]{}

%% use optional labels to link authors explicitly to addresses:
%% \author[label1,label2]{<author name>}
%% \address[label1]{<address>}
%% \address[label2]{<address>}

\author[1]{Evan North\texorpdfstring{\corref{cor1}}{}}
\ead{einorth@ncsu.edu}

\author[1]{Semyon Tsynkov}
\ead{tsynkov@math.ncsu.edu}
\ead[url]{https://stsynkov.math.ncsu.edu}
\cortext[cor1]{Corresponding author}

\author[2]{Eli Turkel}
\ead{turkel@tauex.tau.ac.il}
\ead[url]{http://www.math.tau.ac.il/~turkel/}

\address[1]{Department of Mathematics, North Carolina State University, Box 8205, Raleigh, NC 27695, USA}
\address[2]{School of Mathematical Sciences, Tel Aviv University, Ramat Aviv, Tel Aviv 69978, Israel}

\begin{abstract}
%% Text of abstract
We use the Method of Difference Potentials (MDP) to solve a non-overlapping domain decomposition formulation of the Helmholtz equation. The MDP reduces the Helmholtz equation on each subdomain to a Calderon's boundary equation with projection  on its  boundary. The unknowns for the Calderon's equation are the Dirichlet and Neumann data. Coupling between neighboring subdomains is rendered by applying their respective Calderon's equations to the same data at the common interface. Solutions on individual subdomains are computed concurrently using a straightforward direct solver. We provide numerical examples demonstrating that our method is insensitive to interior cross-points and mixed boundary conditions, as well as large jumps in the wavenumber for transmission problems,  which are known to be problematic for many other Domain Decomposition Methods.
\end{abstract}

\begin{keyword}
Time-harmonic waves   \sep Non-overlapping domain decomposition \sep Calderon's operators \sep Exact coupling between subdomains \sep High-order accuracy \sep Compact finite difference schemes \sep Direct solution \sep Complexity bounds \sep Material interfaces  \sep Interior cross-points
%% keywords here, in the form: keyword \sep keyword

%% MSC codes here, in the form: \MSC code \sep code
%% or \MSC[2008] code \sep code (2000 is the default)

\end{keyword}

\end{frontmatter}

%%
%% Start line numbering here if you want
%%
% \linenumbers

%% main text

\section{Introduction} \label{sec:INTRO}
The Helmholtz equation governs the propagation of time-harmonic waves. For domains many wavelengths in size, it becomes intractable to solve directly, even on modern computers. Non-overlapping Domain Decomposition Methods (DDMs) attempt to alleviate the cost growth by breaking the domain down into smaller, simpler subdomains thus creating subproblems that are coupled to one another along their interfaces. Traditionally, DDMs resolve this coupling by an iterative process that alternates between directly solving a localized approximation of the subproblem and updating the resulting boundary conditions using parameterized transmission conditions. The convergence rate of the iterative process is heavily dependent on the transmission conditions, the choice of which is a highly active research area (see \cite{BOUBENDIR2018, Stolk2013, MATTESI2019, MODAVE2020, gordon2020}, among others). Typically, accounting for the global behavior in the boundary update step --- where the transmission conditions are utilized --- leads to more expensive updates but fewer iterations for convergence. On the other hand, localized approximations in the transmission conditions tend to result in faster updates but more iterations.

The numerical difficulties associated with solving the Helmholtz equation on large domains are further exacerbated if the wavenumber and/or boundary conditions are discontinuous. In this paper, we address the issue of coupling between subdomains by altogether circumventing the iterative process and solving globally for all of the subdomain boundary data. In addition, the wavenumber may undergo jumps across interfaces while the boundary conditions on different segments of the boundary may have different type (mixed boundary conditions). Our methodology yields the exact solution up to the discretization error. The  global behavior  of the solution is accounted for by enforcing the appropriate transmission conditions  at the collections of interfaces (typically, the continuity of the solution itself and its first normal derivative), while solutions on individual subdomains are computed concurrently using a direct solver.

Our approach utilizes several key features of the Method of Difference Potentials (MDP). Originally proposed by Ryaben'kii \cite{ryaben1985boundary,ryaben2002method}, the MDP can be interpreted as a discrete version of the method of Calderon's operators \cite{cald,seel} in the theory of partial differential equations. The MDP reduces a given partial differential equation from its domain to the boundary. The resulting boundary formulation involves an operator equation (Calderon's boundary equation with projection) with both Dirichlet and Neumann data in the capacity of unknowns. Having solved the boundary operator equation, the solution is reconstructed on the domain using Calderon's potential. Therefore, the MDP allows one to parameterize solutions on the domain using their boundary data. This proves very convenient for a domain decomposition framework. Indeed, once the original domain has been partitioned into subdomains, the Calderon's boundary equations for individual subdomains are naturally coupled with the appropriate interface conditions that are also formulated in terms of the Dirichlet and Neumann data. This yields an overall linear system to be solved only at the combined boundary. As a result, the computation of the boundary projection operators is performed ahead of time and completely in parallel.

As the discrete Calderon's operators are pre-computed, the MDP-based domain decomposition appears most convenient to implement in those cases where all subdomains have the same shape. If the wavenumber is also the same everywhere, then the operators are computed only once and subsequently applied to all subdomains. If the wavenumber jumps between subdomains, then the operators are recomputed for each additional value of the wavenumber. The case of identical subdomains implies no limitation of generality though, as the proposed method extends to more elaborate scenarios where the subdomains may differ in shape and the wavenumber may wary inside subdomains. In this paper, however, we focus on congruent subdomains and a piecewise-constant wavenumber. Information on implementing the MDP for general smooth geometries can be found in \cite{Medvinsky2015}. Within the scope of our current setup, we observe the method's performance on Helmholtz transmission problems (particularly those with large jumps in the wavenumber) and domains with cross-points --- points where more than two subdomains meet. Without special consideration, cross-points and jumps in the wavenumber are known to adversely affect the convergence rate of iterative DDMs. While some recent methods have managed to mitigate these effects \cite{BOUBENDIR2018, MODAVE2020, gordon2020}, we emphasize that our method is intrinsically insensitive to these cross-points and jumps in the wavenumber.

The outline of this paper is as follows: In Section~\ref{sec:DDM}, we introduce DDMs for the Helmholtz equation. Section~\ref{sec:MDP} establishes the representative subdomain and covers the necessary information to implement the MDP, ending with the modifications necessary to apply the MDP as a DDM. Details for a practical implementation are outlined in Section~\ref{sec:implementation} and the complexity of the method is discussed in Section~\ref{sec:complexity}. In Section~\ref{sec:numerical}, numerical results are presented to validate the algorithm, corroborate the claims of complexity from Section~\ref{sec:complexity}, and explore the practical limits of the method. Finally, in Section~\ref{sec:conclusions} we provide a summary and propose directions for future research.

%--------------------------------------------%
%
%       Section:
%       Domain Decomposition
%
%--------------------------------------------%

\section{Domain Decomposition} \label{sec:DDM}
\begin{figure}[ht]
    \centering
    \includegraphics[width=0.4\textwidth]{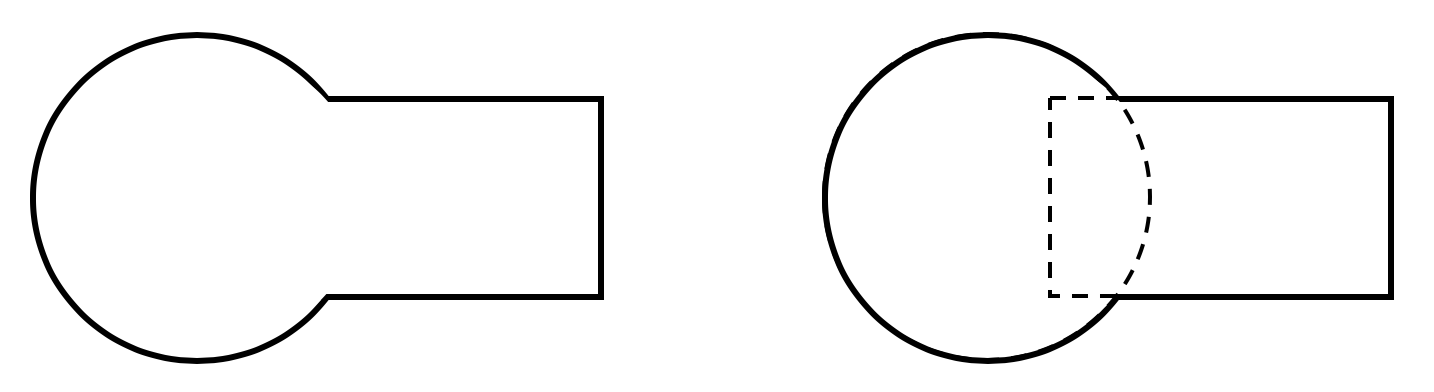}
    \caption{The classical DDM example. An irregular domain composed of two shapes, decomposed into two overlapping subdomains.}
    \label{fig:keyhole}
\end{figure}

Domain Decomposition Methods were originally introduced by Schwarz \cite{schwarz1870} to prove the existence and uniqueness of solutions to the Poisson equation over irregularly shaped domains. The original Schwarz algorithm used overlapping decompositions (Figure~\ref{fig:keyhole}), but was later extended to non-overlapping decompositions (Figure~\ref{fig:simple_domains}) by Lions \cite{Lions1990}. In this paper, we focus on non-overlapping subdomains. Accordingly, we begin with providing a brief overview of non-overlapping DDMs including the original method by Lions for the Poisson equation and subsequent adaptation  by Despr\'es for  the Helmholtz equation. For a more rigorous introduction to DDMs, including proofs of convergence and calculation of convergence factors, see \cite{DOLEAN2015, gander2008}.

%--------------------------------------------%
%
%       Subsection:
%       Non-Overlapping Formulation
%
%--------------------------------------------%

    \subsection{Non-Overlapping Formulation} \label{sec:lions}
\begin{figure}[ht]
        \centering
        \begin{subfigure}{.3\textwidth}
          \centering
            \includegraphics[width=\textwidth] {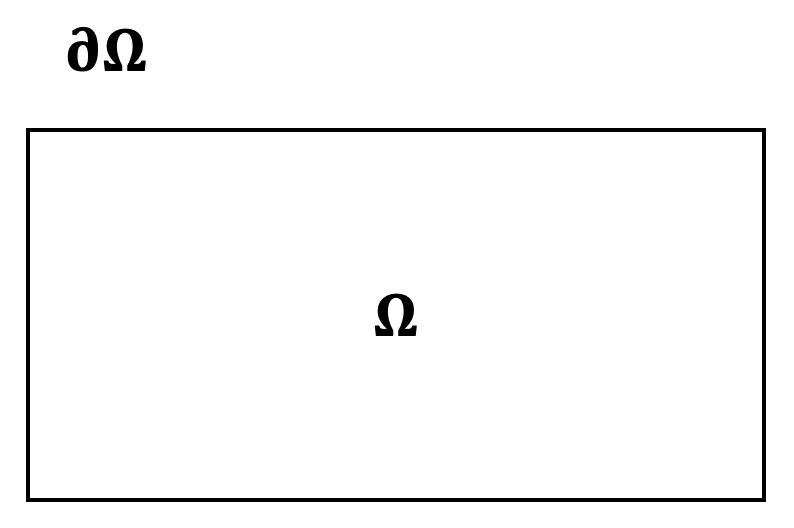}
          \caption{Rectangular domain.}
          \label{fig:simple_rect}
        \end{subfigure}
        \hspace{1.2cm}
        \begin{subfigure}{.3\textwidth}
          \centering
          \includegraphics[width=\textwidth] {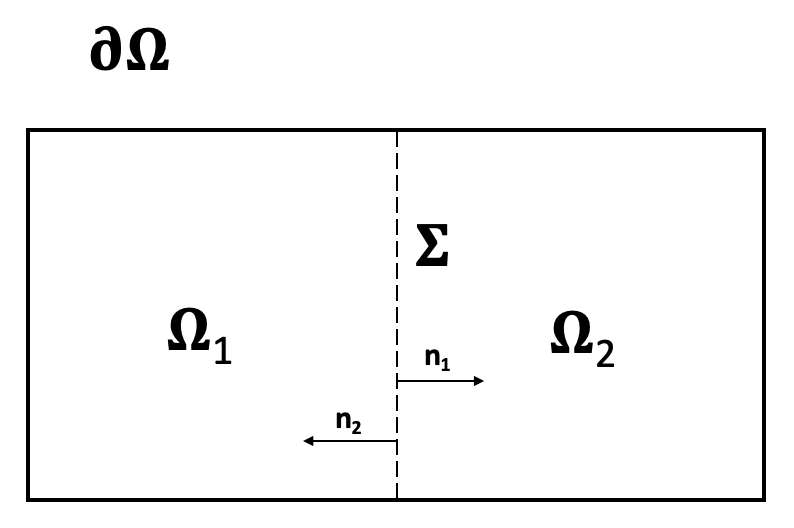}
          \caption{Two-domain decomposition.}
          \label{fig:DDM_rect}
        \end{subfigure}
        \caption{Basic non-overlapping decomposition of a domain $\Omega$ (with boundary $\dO$) into two subdomains, $\Omega_1$ and $\Omega_2$. A fictitious boundary, $\Sigma$, is introduced to indicate the separation between subdomains, and ${\bf n}_i$ is the outward unit normal vector of $\Omega_i$ on $\Sigma$.}
        \label{fig:simple_domains}
    \end{figure}

    Consider the Poisson equation over a rectangular domain $\Omega \subset \R^2$ with boundary $\dO$, (see Figure~\ref{fig:simple_rect}). Then the following Dirichlet boundary value problem (BVP) can be posed:
    \begin{equation}
        \begin{cases} \Delta u = f  \quad &\text{in } \Omega \\
        u = 0 \quad &\texttt{on } \dO \end{cases}
        \label{eq:Poisson_BVP}
    \end{equation}
    Consider a partitioning of $\Omega$ that splits the domain into two subdomains, $\Omega_1$ and $\Omega_2$, by introducing an artificial interface $\Sigma = \overline{\Omega}_1 \cap \overline{\Omega}_2$ as in Figure~\ref{fig:DDM_rect}. The BVP \eqref{eq:Poisson_BVP} can be reformulated over the new subdomains individually:
    \begin{subequations} 
        \label{eq:incomplete_subproblems}
        \begin{align}
         \label{eq:incomplete_subproblems_a}
            &\begin{cases}
                \Delta u_1 = f  & \text{in } \Omega_1 \\
                u_1 = 0 & \text{on } \overline{\Omega}_1 \cap \dO \\
            \end{cases} \\
             \label{eq:incomplete_subproblems_b}
            &\begin{cases}
                \Delta u_2 = f  & \text{in } \Omega_2 \\
                u_2 = 0 & \text{on } \overline{\Omega}_2 \cap \dO \\
            \end{cases}\\
            \label{eq:incomplete_subproblems_c}
            &\begin{cases}
                u_1 = u_2, & \text{on } \Sigma \\
                \partialdiff{u_1}{\bf{n_1}} = - \partialdiff{u_2}{\bf{n_2}}& \text{on } \Sigma \\
            \end{cases}
        \end{align}
    \end{subequations}
    where the interface conditions (\ref{eq:incomplete_subproblems_c}) guarantee that the combined solution of (\ref{eq:incomplete_subproblems})  coincides with that of \eqref{eq:Poisson_BVP}:
    \begin{equation*}
        \begin{cases}
            u_1 = u, & \text{in } \Omega_1 \\
            u_2 = u, & \text{in } \Omega_2 \\
        \end{cases}
    \end{equation*}
    Conditions other than (\ref{eq:incomplete_subproblems_c}) can be formulated on $\Sigma$ so that the resulting combined problem is well-posed, but its solution will be different from the true solution of \eqref{eq:Poisson_BVP}.

    There are two separate interface conditions in (\ref{eq:incomplete_subproblems_c}). They apply to both subproblems (\ref{eq:incomplete_subproblems_a}) and (\ref{eq:incomplete_subproblems_b}) at the same time and couple these subproblems together. However, each subproblem (\ref{eq:incomplete_subproblems_a}) or (\ref{eq:incomplete_subproblems_b}) considered independently, i.e., with no connection to the other subproblem, is not fully specified and cannot be solved on its own because it is missing boundary conditions on $\Sigma$. To enable the individual solvability, one needs to provide these boundary conditions. Yet unlike in (\ref{eq:incomplete_subproblems_c}), one cannot specify more than one boundary condition on $\Sigma$ for either of the two standalone problems (\ref{eq:incomplete_subproblems_a}) or (\ref{eq:incomplete_subproblems_b}), as that would result in an overdetermination. In other words, when solving (\ref{eq:incomplete_subproblems_a}) one cannot specify both $u_1$ and $\partialdiff{u_1}{\bf{n_1}}$ on $\Sigma$, and likewise for  (\ref{eq:incomplete_subproblems_b}).% and $u_2$.

    To avoid the overdetermination 
    and still allow for separate solution of individual subproblems, 
    P.L. Lions proposed to use one Robin boundary condition \cite{Lions1990},  formed as a linear combination of the two continuity conditions (\ref{eq:incomplete_subproblems_c}). For any pair of constants $\left( p_1, p_2 \right) \in \R^2$, this \textit{transmission condition} yields the following combined formulation in lieu of (\ref{eq:incomplete_subproblems}):
    \begin{subequations}
        \label{eq:Lions_subproblems}
        \begin{align}
            \label{eq:Lions_subproblems_a}
            &\begin{cases}
                \Delta u_1 = f  \quad \text{in } \Omega_1 \\
                u_1 = 0 \quad \text{on } \overline{\Omega}_1 \cap \dO \\
                \left( \partialdiff{}{\bf n_1} + p_1 \right) u_1 = \left( \partialdiff{}{\bf n_1} + p_1 \right) u_2 \quad \text{on } \Sigma
            \end{cases} \\
            \label{eq:Lions_subproblems_b}
            &\begin{cases}
                \Delta u_2 = f  \quad \text{in } \Omega_2 \\
                u_2 = 0 \quad \text{on } \overline{\Omega}_2 \cap \dO \\
                \left( \partialdiff{}{\bf n_2} + p_2 \right) u_2 = \left( \partialdiff{}{\bf n_2} + p_2 \right) u_1 \quad \text{on } \Sigma
            \end{cases}
        \end{align}
    \end{subequations}
    Each of the two subproblems (\ref{eq:Lions_subproblems}) is individually well-defined in the sense that the third equation in either (\ref{eq:Lions_subproblems_a}) or (\ref{eq:Lions_subproblems_b}) can be interpreted as a Robin boundary condition on $\Sigma$ for $u_1$ or $u_2$, respectively, with the right-hand side of the respective equation providing the data. However, the relation of the combined formulation (\ref{eq:Lions_subproblems}) to  the original BVP \eqref{eq:Poisson_BVP} requires a special inquiry.
    
    Lions conducted the corresponding analysis in \cite{Lions1990}.  He replaced  the combined formulation (\ref{eq:Lions_subproblems}) with the iteration:
    \begin{align}
        \begin{cases} \Delta u_1^{(n+1)} = f  \quad \text{in } \Omega_1 \\
        u_1^{(n+1)} = 0 \quad \text{on } \overline{\Omega}_1 \cap \dO \\
        \left( \partialdiff{}{\bf n_1} + p_1 \right) u_1^{(n+1)} = \left( \partialdiff{}{\bf n_1} + p_1 \right) u_2^{(n)} \quad \text{on } \Sigma
        \end{cases} \\
        \begin{cases} \Delta u_2^{(n+1)} = f  \quad \text{in } \Omega_2 \\
        u_2^{(n+1)} = 0 \quad \text{on } \overline{\Omega}_2 \cap \dO \\
        \left( \partialdiff{}{\bf n_2} + p_2 \right) u_2^{(n+1)} = \left( \partialdiff{}{\bf n_2} + p_2 \right) u_1^{(n)} \quad \text{on } \Sigma
        \end{cases}
        \label{eq:two_subproblems_iterative}
    \end{align}
    and proved that this iteration converges to the solution of \eqref{eq:Poisson_BVP} as $n$ increases. The rate of convergence depends on the choice of the parameters $p_1$ and $p_2$. As the next iteration $n+1$ for each subproblem only relies on the other subproblem's current iteration $n$, the subproblems can be solved in parallel to one another, a highly desirable trait for DDMs. The proof given in \cite{Lions1990} extends to an arbitrary number of subdomains.

%--------------------------------------------%
%
%       Subsection:
%       Helmholtz Adaptation
%
%--------------------------------------------%

    \subsection{Helmholtz Adaptation} \label{sec:despres}
    Complications arise when applying \eqref{eq:two_subproblems_iterative} directly to the Helmholtz equation. Consider the following BVP over the domain from Figure~\ref{fig:simple_rect}:
    \begin{equation}
        \begin{cases} \Delta u + k^2 u = f  \quad &\text{in } \Omega \\
        u = 0 \quad &\texttt{on } \dO \end{cases}
        \label{eq:helmholtz_BVP}
    \end{equation}
    To guarantee well-posedness of \eqref{eq:helmholtz_BVP}, i.e., to avoid resonance, $-k^2$ may not be an eigenvalue of the underlying Laplace problem. However, when considering a decomposition such as the one in Figure~\ref{fig:DDM_rect} with the Lions transmission condition, it is non-trivial to know that $-k^2$ will always remain outside the spectrum of the corresponding Laplace subproblem, which only becomes more problematic when various decompositions are  considered. This issue was addressed in \cite{despres_1993} when Despr\'es proposed the use of Lions' transmission condition with $p_1 = p_2 = ik$ (where $i = \sqrt{-1}$). This choice yields the following subproblems (cf.\ (\ref{eq:Lions_subproblems})):
    \begin{subequations}
        \label{eq:two_subproblems_despres}
        \begin{align}
                    \label{eq:two_subproblems_despres_a}
            \begin{cases} (\Delta + k^2) u_1 = f  \quad \text{in } \Omega_1 \\
            u_1 = 0 \quad \text{on } \overline{\Omega}_1 \cap \dO \\
            \left( \partialdiff{}{\bf n_1} + ik \right) u_1 = \left( \partialdiff{}{\bf n_1} + ik \right) u_2 \quad \text{on } \Sigma
            \end{cases} \\
            \label{eq:two_subproblems_despres_b}
            \begin{cases} (\Delta + k^2) u_2 = f  \quad \text{in } \Omega_2 \\
            u_2 = 0 \quad \text{on } \overline{\Omega}_2 \cap \dO \\
            \left( \partialdiff{}{\bf n_2} + ik \right) u_2 = \left( \partialdiff{}{\bf n_2} + ik \right) u_1 \quad \text{on } \Sigma
            \end{cases}
        \end{align}
    \end{subequations}
    Despr\'es' transmission condition  shifts the spectrum of the operator to the complex domain, guaranteeing that resonant frequencies are avoided on each subproblem (\ref{eq:two_subproblems_despres_a}) or (\ref{eq:two_subproblems_despres_b}). It does so at the cost of introducing complex values into the problem, but for many applications, this is computationally not an issue. An iterative procedure similar to (\ref{eq:two_subproblems_iterative}) can be employed for \eqref{eq:two_subproblems_despres}, and Despr\'es showed in \cite{despres_1993} that it will converge.

%--------------------------------------------%
%
%       Subsection:
%       Other Considerations
%
%--------------------------------------------%

    \subsection{Other Considerations}
    The methods outlined above are the foundation of most modern DDMs for the Helmholtz equation, and have been improved upon in recent years. For example, quasi-optimal convergence rates have been achieved by optimizing the choice of transmission conditions with the so-called ``square root operator" \cite{BOUBENDIR2012}. However, while this leads to convergence in fewer iterations, it generally requires more expensive iterations.

    Recent work has also been dedicated to the resolution of interior cross-points. The interior cross-points are points where more than two subdomains meet, and they pose no issues at the continuous level of the formulation. Yet the cross-points are known to adversely affect the accuracy and convergence if not discretized with care. In \cite{gander2016}, several methods are discussed for resolving these cross-points for elliptic problems, and \cite{MODAVE2020} provides an extension of the quasi-optimal method from \cite{BOUBENDIR2012} that accounts for interior cross-points. In Section~\ref{sec:numerical}, we demonstrate how the issue of cross-points is resolved naturally with our method, with no special consideration. While some other methods can also address the cross-points (see, for example, \cite{gordon2020}), we emphasize that our method is completely insensitive to them by design.

    Transmission problems provide another common venue for the application of DDMs, but they can  require special care in the high-contrast, high-frequency regime (see \cite{BOUBENDIR2018}, an extension of the square-root operator from \cite{BOUBENDIR2012}). Similarly to the case of cross-points, our method appears to be insensitive to large jumps in the wavenumber, as discussed further in Section~\ref{sec:piecewise_k}.

%--------------------------------------------%
%
%       Section:
%       Method of Difference Potentials
%
%--------------------------------------------%

\section{Method of Difference Potentials} \label{sec:MDP}
To introduce the Method of Difference Potentials \cite{ryaben2002method}, consider the inhomogeneous Helm\-holtz equation with a general (constant-coefficient) Robin boundary condition
\begin{subequations} \label{eq:helmholtz_physical}
    \begin{equation} \label{eq:helmholtz}
        \Delta u + k^2 u = f
    \end{equation}
    \begin{equation} \label{eq:helmholtz_robin}
        \alpha u + \beta \partialdiff{u}{\bf n} = \phi
    \end{equation}
\end{subequations}
over the domain $\Omega \subset \R^2$ depicted in Figure~\ref{fig:simple_rect}, as well as its decomposition depicted in Figure~\ref{fig:DDM_rect}. In a similar manner to traditional DDMs, we split the problem into two separate subdomains as in \eqref{eq:incomplete_subproblems}, and encounter the same issue of needing to enforce continuity of the solution and its flux over the interface $\Sigma$.

The key role of the MDP is to impose the required interface conditions on $\Sigma$. The MDP replaces the governing differential equation, the Helmholtz equation \eqref{eq:helmholtz}, on the domain by an equivalent operator equation at the boundary (Calderon's boundary equation with projection). The latter is formulated with respect to the Cauchy data of the solution, i.e., the boundary trace of the solution itself (Dirichlet data) and its normal derivative (Neumann data). The reduction to the boundary is done independently for individual subdomains $\Omega_1$ and $\Omega_2$ (see Figure~\ref{fig:DDM_rect}). Then, the resulting boundary equations with projections on the neighboring subdomains share the Dirichlet and Neumann data as unknowns at the common interface $\Sigma$, which directly enforces the continuity of the solution and its flux. For the remaining parts of the boundaries, $\partial\Omega_1\backslash\Sigma$ and $\partial\Omega_2\backslash\Sigma$, the boundary equations with projections are combined with the boundary condition \eqref{eq:helmholtz_robin}, which is  formulated in terms of the Cauchy data of the solution. Altogether, the MDP solves a fully coupled problem for the Helmholtz equation  similar to \eqref{eq:incomplete_subproblems}. Nonetheless,  it turns out that the constituent subproblems can still be solved independently and in parallel. 

\begin{figure}[ht]
    \centering
    \includegraphics[width=0.7\textwidth]{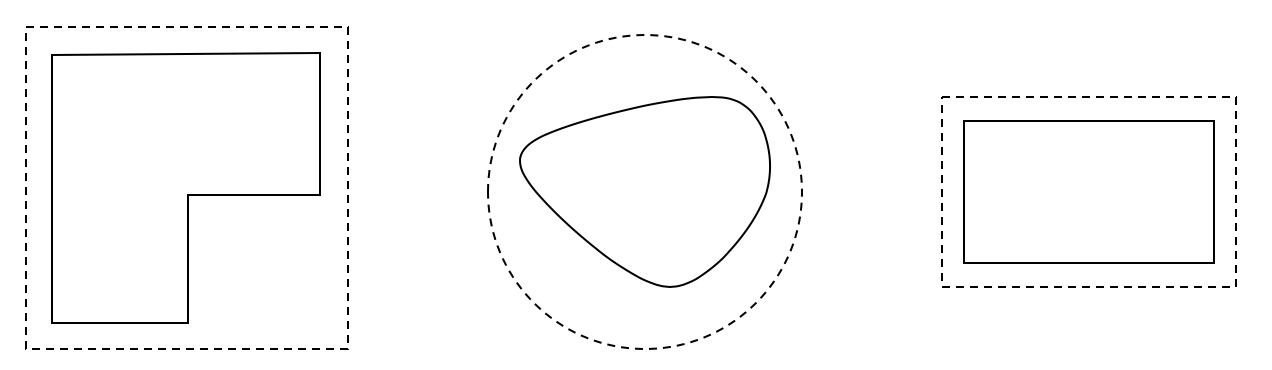}
    \caption{Example domains (solid border) and a reasonable choice of auxiliary domain (dotted border) for each.}
    \label{fig:aux_examples}
\end{figure}

\begin{figure}[ht]
    \centering
    \includegraphics[width=0.7\textwidth]{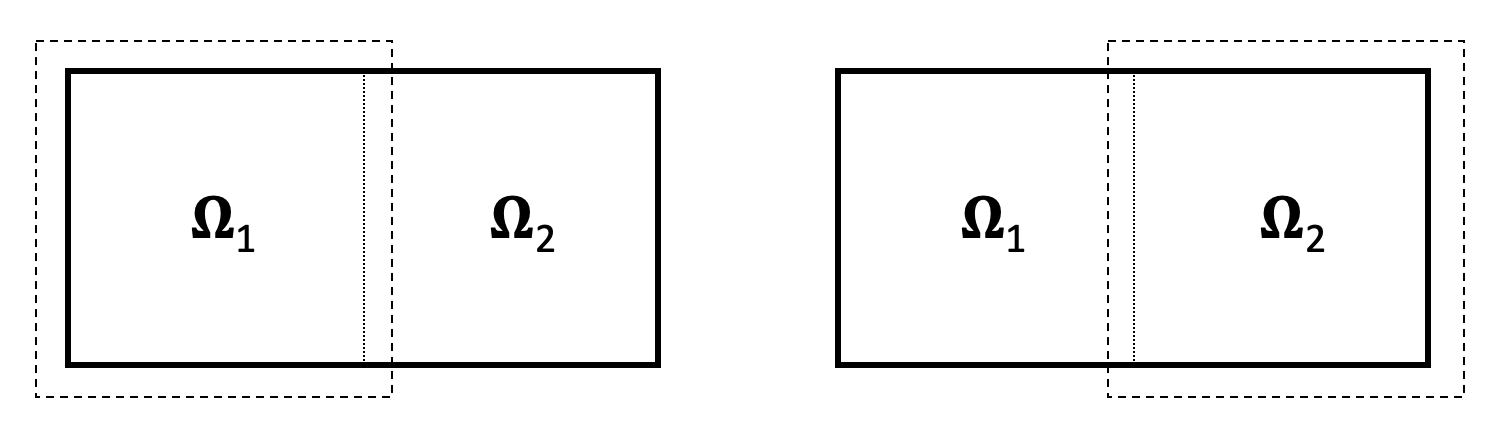}
    \caption{The auxiliary domain setup for our problem with the domain decomposition from Figure~\ref{fig:DDM_rect}.}
    \label{fig:aux_ddm}
\end{figure}

Indeed, Calderon's operators are computed with the help of the auxiliary problem, which is formulated for the same governing equation, but on a larger auxiliary domain. The auxiliary problem must be uniquely solvable and well-posed. Otherwise, the auxiliary problem can be arbitrary, and is normally chosen so as to enable an easy and efficient numerical solution. In particular, the auxiliary domain would typically have a simple regular shape; some examples are shown in Figure~\ref{fig:aux_examples}. Given that for domain decomposition one needs to compute the Calderon operators separately for individual subdomains, we embed each subdomain within its own auxiliary domain, see Figure~\ref{fig:aux_ddm}, and solve the resulting auxiliary problems independently. In practice, we also take advantage of the fact that in some of our simulations the subdomains are identical and reuse the computed operators accordingly. Precise criteria for the selection of an auxiliary domain, as well as details of how to efficiently account for identical subdomains, are discussed in Section~\ref{sec:auxiliary_problem}.

In the rest of this section, we introduce the parts of the MDP necessary to implement it in the framework of DDM. For a detailed account of the theory and derivation of the MDP, see \cite{ryaben2002method}, as well as \cite{Medvinsky2012, Britt_10}, among others. Additionally, for details on handling more complicated boundary conditions, as well as extending this method to domains with curvilinear sides, see \cite{Britt2013} or \cite{Medvinsky2015}, respectively.

%--------------------------------------------%
%
%       Subsection:
%       Finite Difference Scheme
%
%--------------------------------------------%
    \subsection{Finite Difference Scheme} \label{sec:FDM}
\begin{figure}[ht]
        \centering
        \includegraphics[width=0.4 \textwidth]{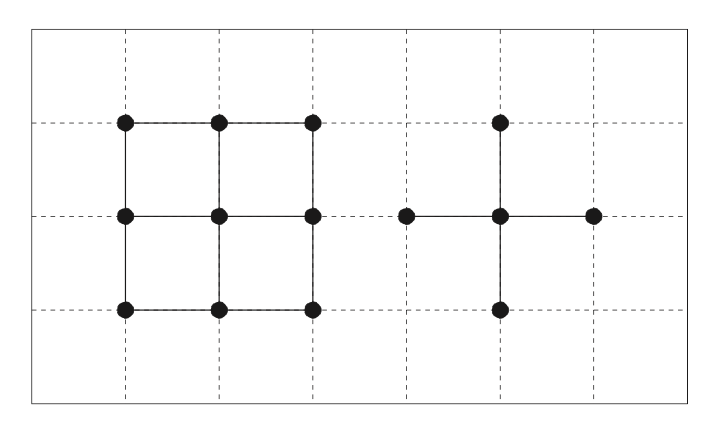}
        \caption{The stencils for the compact scheme given in \eqref{eq:FD_scheme}.}
        \label{fig:FD_scheme}
    \end{figure}

    The MDP can be implemented in conjunction with any finite difference scheme as the underlying approximation, including the case of complex or non-conforming boundaries \cite{Medvinsky2015}. High-order schemes are known to reduce the pollution effect for the Helmholtz equation \cite{Babuska_00, Bayliss_83, Deraemaker_99}. Further, compact schemes require no additional boundary conditions beyond what is needed for the differential equation itself. Therefore, we have chosen to use the fourth-order, compact scheme for the Helmholtz equation as presented in \cite{harari_turkel_1995, singer_turkel_1998}:
    \begin{align} \label{eq:FD_scheme}
        &\frac{1}{h^2} \left( u_{m+1,n} + u_{m-1,n} + u_{m,n+1} + u_{m,n-1} - 4u_{m,n} \right) \nonumber \\
        &+ \frac{1}{6h^2} \left[ u_{m+1,n+1} + u_{m-1,n+1} + u_{m+1,n-1} + u_{m-1,n-1} + 4u_{m,n} \nonumber \right. \\
        &- \left. 2\left( u_{m+1,n} + u_{m-1,n} + u_{m,n+1} + u_{m,n-1} \right) \right] \\
        &+ \frac{k^2}{12} \left( u_{m+1,n} + u_{m-1,n} + u_{m,n+1} + u_{m,n-1} + 8u_{m,n} \right) \nonumber \\
        &= f_{m,n} + \frac{1}{12} \left( f_{m+1,n} + f_{m-1,n} + f_{m,n+1} + f_{m,n-1} - 4f_{m,n} \right) \nonumber
    \end{align}
    The scheme in \eqref{eq:FD_scheme} uses a nine-node stencil for the left-hand side of the PDE and a five-node stencil for the right-hand side (see Figure~\ref{fig:FD_scheme}), with uniform step size in both directions ($\Delta x = \Delta y = h$). In the case where the PDE is homogeneous, the right-hand side stencil is unnecessary as $f \equiv 0$. Additionally, \eqref{eq:FD_scheme} was derived for a constant value of the wavenumber $k$. For our purposes this is sufficient, because while the domain $\Omega$ may have a piecewise-constant $k$, we assume that the decomposition is such that each $\Omega_i$ has constant $k$. One could also consider a sixth-order scheme for constant \cite{Turkel_13} or variable \cite{Singer_11} wavenumber $k$, or a fourth-order scheme for a more general form of the Helmholtz equation with a variable coefficient Laplace-like term and wavenumber \cite{Britt_10}. However, for the scope of this paper we will focus on the piecewise constant $k$ case.

%--------------------------------------------%
%
%       Subsection:
%       Base Subdomain
%
%--------------------------------------------%

    \subsection{Base Subdomain} \label{sec:base_subdomain}
    In this paper, we focus on situations where the problem can be decomposed into identical subdomains so that as much information as possible can be reused. We define all of the components needed to perform the MDP algorithm on one base subdomain, and allow copies of that base subdomain to be translated and rotated into the appropriate position for any given concrete example. Considering the model domain from Figure~\ref{fig:simple_domains}, a logical choice of base subdomain is a square. For the sake of introducing the MDP on the base subdomain, throughout Section~\ref{sec:base_subdomain} we will refer to the base subdomain simply as $\Omega$ with boundary $\Gamma$, where $\Omega$ is a square with a side length of 2, centered at the origin. This simple cased is used for efficiency. There is no substantial difficulty to treat subdivisions that are different or even have non-rectangular shape.

%--------------------------------------------%
%
%       Subsubsection:
%       Auxiliary Problem
%
%--------------------------------------------%
        \subsubsection{Auxiliary Problem} \label{sec:auxiliary_problem}
        We will embed the base subdomain $\Omega$ in a larger domain $\Omega_0$. This larger domain is known as the auxiliary domain, on which we formulate  the auxiliary problem (AP). The AP should be uniquely solvable and well-posed, and should  allow for a  convenient and efficient numerical solution.

        Let $\L$ represent the Helmholtz operator: $\L u \defeq (\Delta + k^2) u$. We formulate the AP on $\Omega_0$ by supplementing the inhomogeneous Helmholtz equation  with homogeneous Dirichlet conditions on the $y-$boundaries and local Sommerfeld conditions on the $x-$boundaries:
        \begin{equation}
            \begin{cases}
                \L u = g,  &(x,y) \in \Omega_0 \\
                u = 0, &y = \pm 1.1 \\
                \frac{\partial u}{\partial x} + iku = 0, &x = 1.1 \\
                \frac{\partial u}{\partial x} - iku = 0, &x = -1.1
            \end{cases}
            \label{eq:auxiliary_problem}
        \end{equation}
        The choice of Sommerfeld-type conditions on the $x-$boundaries makes the spectrum of the AP (\ref{eq:auxiliary_problem}) complex, guaranteeing that resonance is avoided for any real wavenumbers $k$. Hence, the AP (\ref{eq:auxiliary_problem}) has a unique solution $u$ for any right-hand side $g$.  It should be noted that although similar in form to the Despr\'es condition from Section~\ref{sec:despres}, the Sommerfeld-type conditions in (\ref{eq:auxiliary_problem}) do not serve any transmission-related purpose, as they exist solely on the auxiliary domain and not on the physical boundary $\Gamma$.

        To discretize the  AP (\ref{eq:auxiliary_problem}), we first replace the operator $\L$ with the left-hand side of the scheme  \eqref{eq:FD_scheme}:
        \begin{subequations} \label{eq:discrete_AP}
            \begin{align} \label{eq:aux_FD_scheme}
                &\frac{1}{h^2} \left( u_{m+1,n} + u_{m-1,n} + u_{m,n+1} + u_{m,n-1} - 4u_{m,n} \right) \nonumber \\
                &+ \frac{1}{6h^2} \left[ u_{m+1,n+1} + u_{m-1,n+1} + u_{m+1,n-1} + u_{m-1,n-1} + 4u_{m,n}  \right. \\
                &- \left. 2\left( u_{m+1,n} + u_{m-1,n} + u_{m,n+1} + u_{m,n-1} \right) \right] \nonumber \\
                &+ \frac{k^2}{12} \left( u_{m+1,n} + u_{m-1,n} + u_{m,n+1} + u_{m,n-1} + 8u_{m,n} \right) = g_{m,n} . \nonumber
            \end{align}
            To maintain the overall accuracy of the solution, the boundary conditions also need to be approximated to fourth-order. For the $y-$boundaries this is trivial, as the boundary nodes can directly be set to zero, i.e. for $m = 0,\, ...,\, M$ set
            \begin{equation} \label{eq:AP_bounds_top_bottom}
                u_{m,0} = u_{m,N} = 0
            \end{equation}
            The following discretization of the Sommerfeld-type conditions was derived for the variable coefficient Helmholtz equation in \cite{Britt_10} and simplified for the constant coefficient case in \cite{Britt2013}:
            \begin{align}
                \begin{split} \label{eq:AP_bounds_right}
                    &\left( \frac{u_{M,n} - u_{M-1,n}}{h} - \frac{1}{6h} \left( u_{M,n+1} - u_{M-1,n+1} + u_{M,n-1} - u_{M-1, n-1} \right. \right.\\
                    &\left.\left. -2 \left( u_{M,n} - u_{M-1,n}\right) \right) - \frac{k^2 h}{24} \left( u_{M,n} - u_{M-1,n} \right) \right) \\
                    &+ik \left( \frac{u_{M,n} - u_{M-1,n}}{h} + \frac{h^2 k^2}{8} u_{M-\frac{1}{2},n} \right.\\
                    &\left. + \frac{u_{M-\frac{1}{2}, n+1} - 2u_{M-\frac{1}{2}, n} + u_{M-\frac{1}{2}, n-1}}{2} \right) = 0
                \end{split}\\
                \begin{split} \label{eq:AP_bounds_left}
                    &\left( \frac{u_{1,n} - u_{0,n}}{h} - \frac{1}{6h} \left( u_{1,n+1} - u_{0,n+1} + u_{1,n-1} - u_{0, n-1} -2 \left( u_{1,n} - u_{0,n}\right) \right) \right.\\
                    &\left.  - \frac{k^2 h}{24} \left( u_{1,n} - u_{0,n} \right) \right) \\
                    &-ik \left( \frac{u_{1,n} - u_{0,n}}{h} + \frac{h^2 k^2}{8} u_{\frac{1}{2},n} +\frac{u_{\frac{1}{2}, n+1} - 2u_{\frac{1}{2}, n} + u_{\frac{1}{2}, n-1}}{2} \right) = 0 .
                \end{split}
            \end{align}
        \end{subequations}
        Conditions \eqref{eq:AP_bounds_right} and \eqref{eq:AP_bounds_left} were derived under the assumption that the source function is compactly supported. In our current setting, the  grid function $g_{m,n}$ will be specified on the interior grid nodes, $m = 1,\, ...,\, M-1$ and $n = 1,\, ...,\, N-1$, and will be zero on the outermost grid nodes.

        We define the discrete operator $\Lh$ as the application of the left-hand side of \eqref{eq:aux_FD_scheme}, allowing the discrete AP to be expressed as $\Lh u = g$ subject to the boundary conditions from \eqref{eq:AP_bounds_top_bottom}, \eqref{eq:AP_bounds_right}, and \eqref{eq:AP_bounds_left}. Similar to the continuous AP \eqref{eq:auxiliary_problem}, the  finite difference AP \eqref{eq:discrete_AP} has a unique solution $u$ for any discrete right-hand side $g$. This solution $u$ defines the inverse operator $\Gh$: $u = \Gh g$.

        In particular, the right-hand side  $g$ may be defined  as
        \begin{equation}
            \label{eq:rhs_g}
            g_{m,n} = \Bh f_{m,n} \defeq f_{m,n} + \frac{1}{12} \left( f_{m+1,n} + f_{m-1,n} + f_{m,n+1} + f_{m,n-1} - 4f_{m,n} \right)
        \end{equation}
        where $\Bh$ represents the application of the stencil from the right-hand side of the scheme \eqref{eq:FD_scheme}. We emphasize that $\Gh$ is defined for any grid function $g$, not just those of the form $g = \Bh f$. The discrete AP can be solved by a combination of a sine-FFT in the $y-$direction and a tridiagonal solver in the $x-$direction to create an efficient approximation method for the solution to the continuous AP \eqref{eq:auxiliary_problem}.

%--------------------------------------------%
%
%       Subsubsection:
%       Grid Sets and Difference Potentials
%
%--------------------------------------------%

        \subsubsection{Grid Sets and Difference Potentials} \label{sec:grid_sets_and_operators}
\begin{figure}[H]
            \centering
            \begin{subfigure}{0.3\textwidth}
                \centering \includegraphics[width=\textwidth]{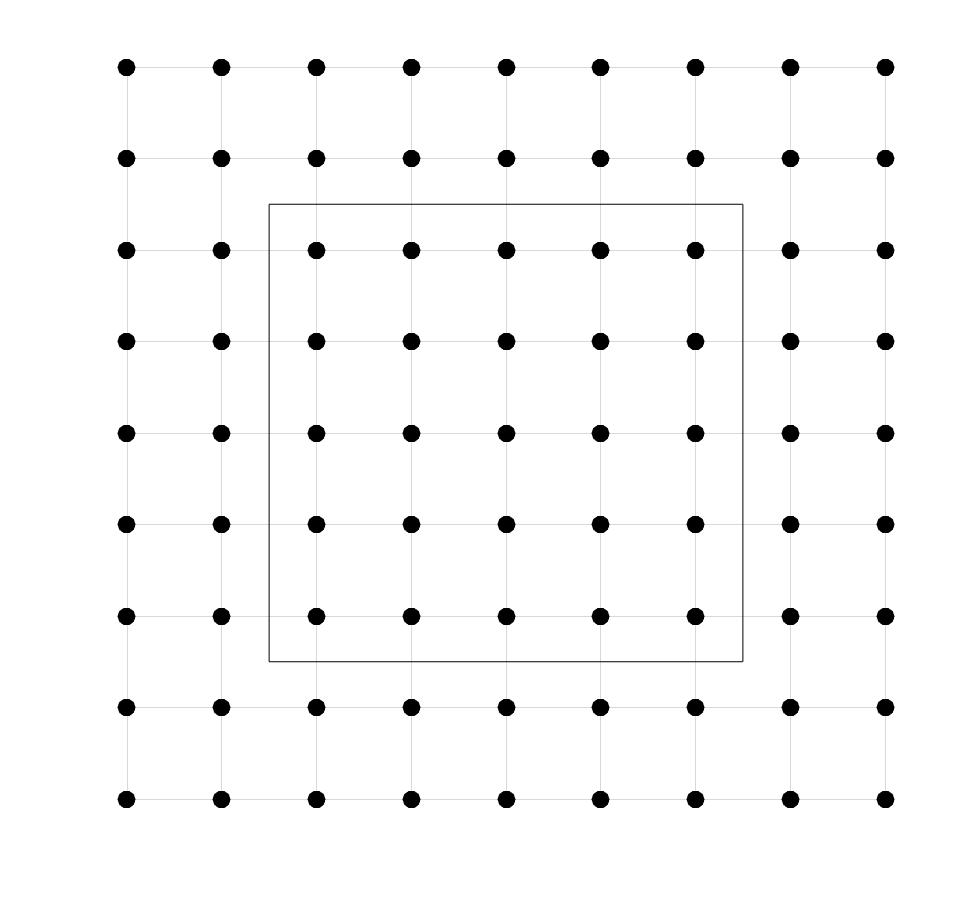}
                \caption{$\N_0$}
                \label{fig:N0}
            \end{subfigure} \hspace{1cm}
            \begin{subfigure}{0.3\textwidth}
                \centering
                 \includegraphics[width=\textwidth]{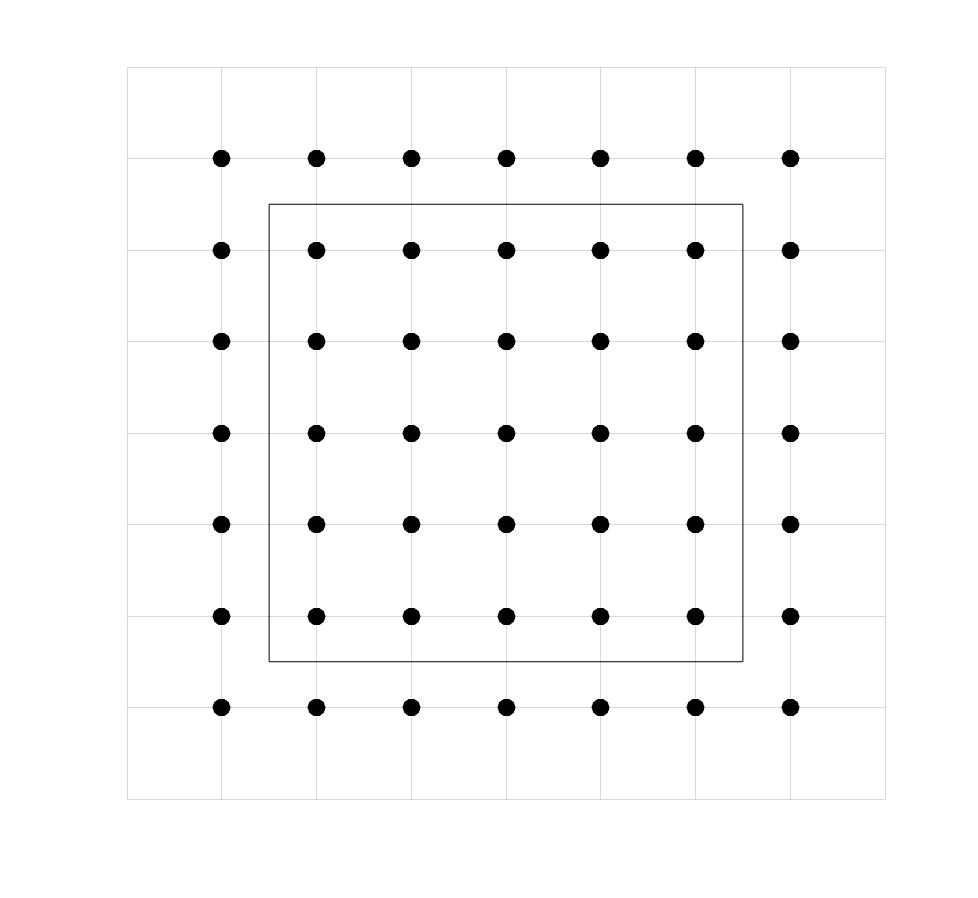}
                 \caption{$\M_0$}
                 \label{fig:M0}
            \end{subfigure}
            \caption{Cartesian grid sets used for the stencils presented in Figure~\ref{fig:FD_scheme} overlaid with the domain $\Omega$.}
            \label{fig:M0_N0}
        \end{figure}

\begin{figure}[H]
        \centering
            \begin{subfigure}{.3\textwidth}
                \includegraphics[width=\textwidth]{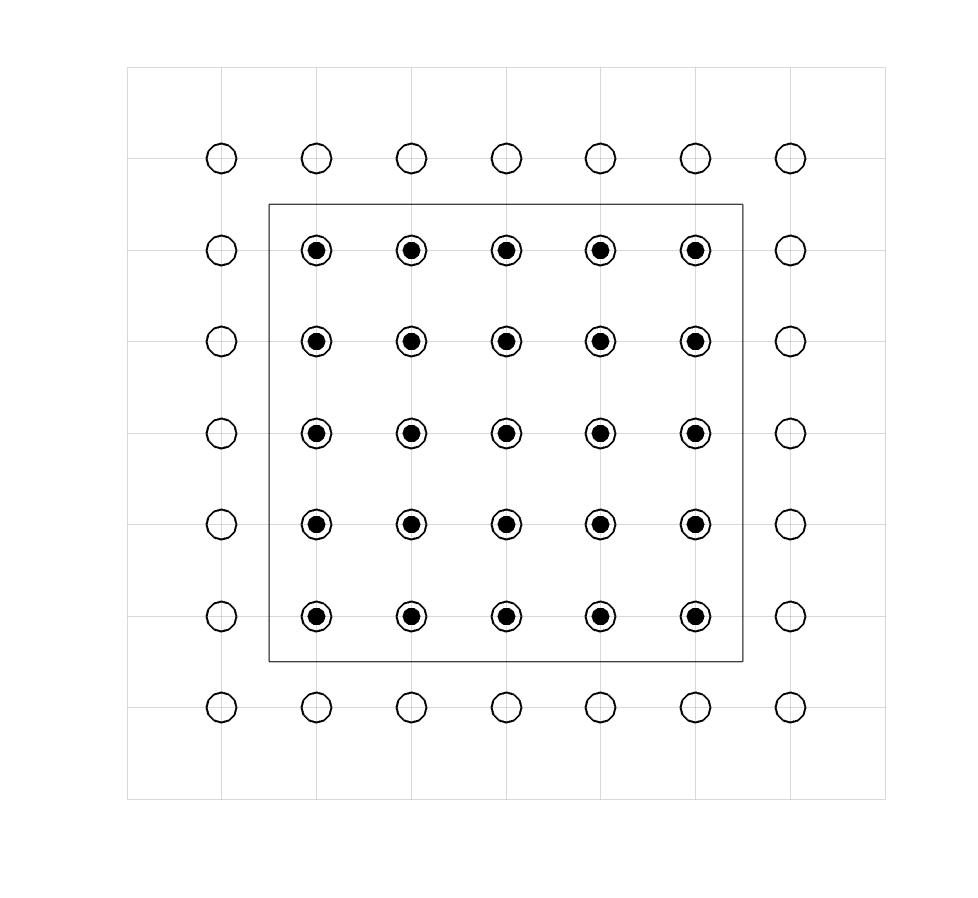}
                \caption{$\cdot-\M^+$, $\circ-\N^+$}
                \label{fig:np_mp}
            \end{subfigure}
            \begin{subfigure}{.3\textwidth}
                \includegraphics[width=\textwidth]{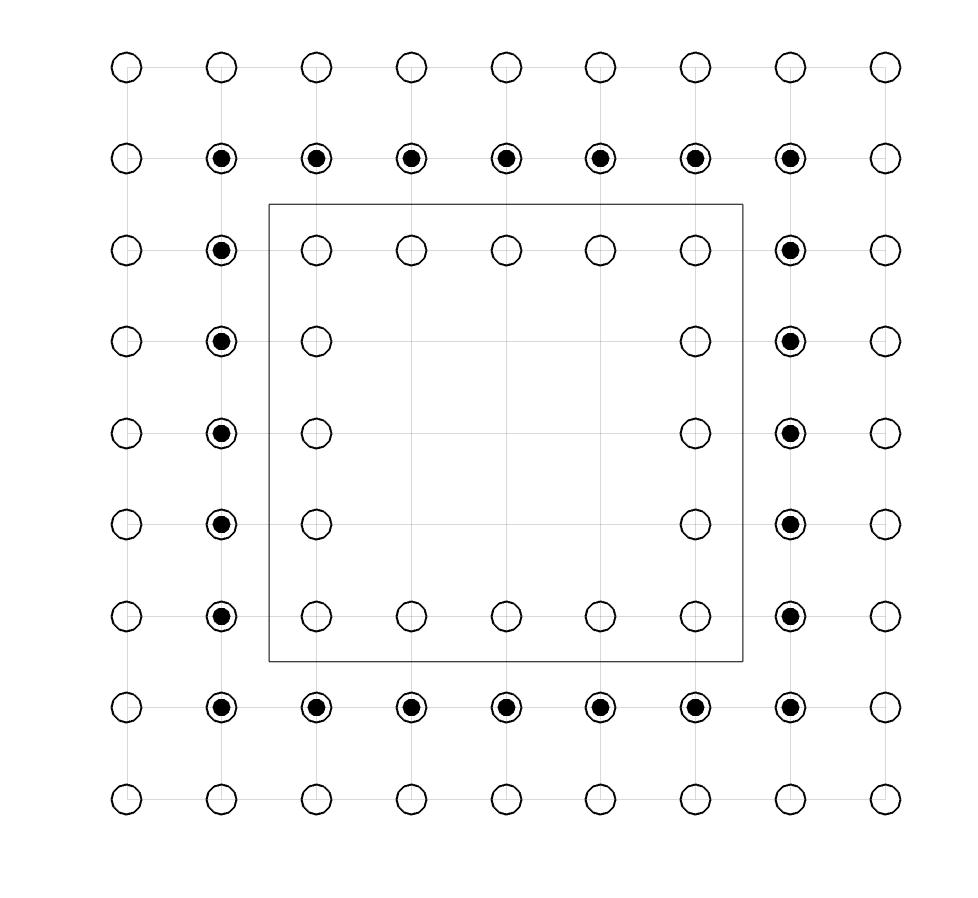}
                \caption{$\cdot-\M^-$, $\circ-\N^-$}
                \label{fig:nm_mm}
            \end{subfigure}
            \begin{subfigure}{.3\textwidth}
                \includegraphics[width=\textwidth]{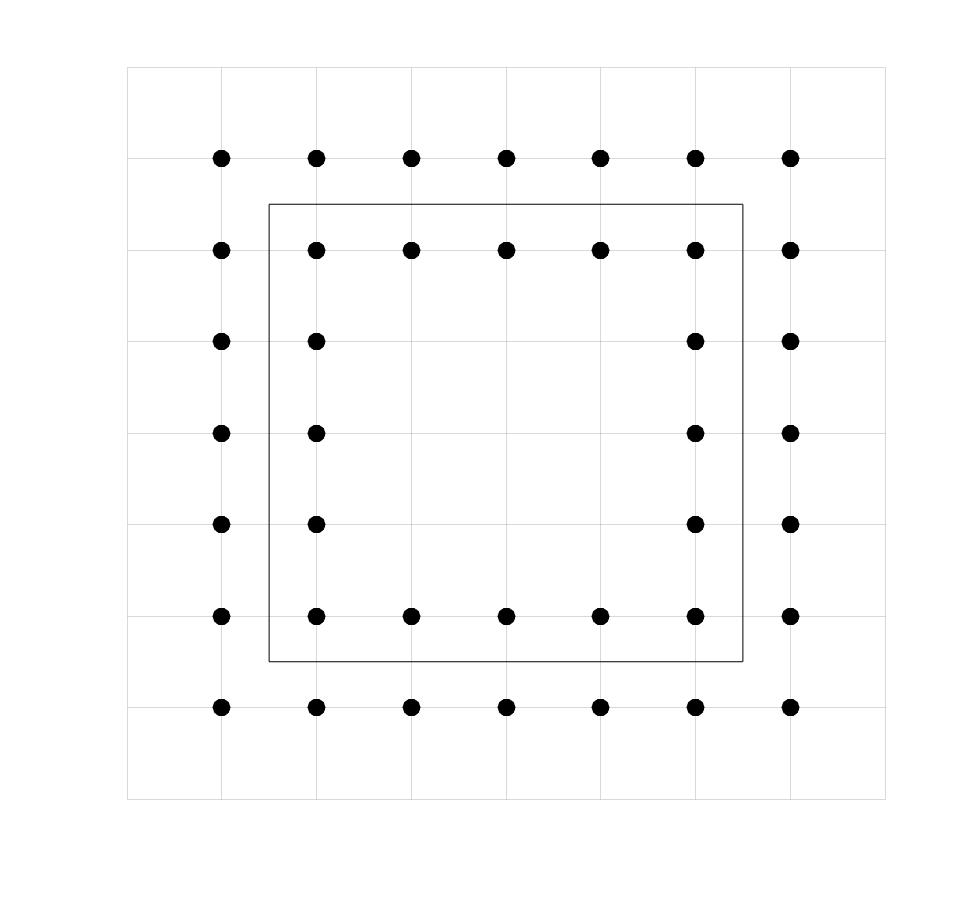}
                \caption{$\cdot-\gamma$}
                \label{fig:gamma}
            \end{subfigure}
            \caption{Discrete analogue grid sets of the interior, exterior, and boundary with respect to the original problem domain $\Omega$.}
            \label{fig:grid_sets}
        \end{figure}

        Let $\N_0$ be a  Cartesian grid  on $\Omega_0$ with uniform step size $h$ in both the $x-$ and $y-$ directions. Let $\M_0 \subset \N_0$ be the set of nodes strictly interior to $\Omega_0$, i.e. not on the boundary (see Figure~\ref{fig:M0_N0}). Define $\M^+ = \M_0 \cap \Omega$ as the nodes that are interior to the original domain $\Omega$, and the exterior nodes as $\M^- = \M_0 \backslash \M^+$. Let $\N^+$ be the set of nodes needed to apply the $3\times 3$ stencil from Figure~\ref{fig:FD_scheme}(left) to every node in $\M^+$, and similarly let $\N^-$ be the same for $\M^-$ (see Figures~\ref{fig:np_mp} and \ref{fig:nm_mm}). Finally, we define the grid boundary  $\gamma = \N^+ \cap \N^-$ as the discrete analogue of the original problem's boundary, $\Gamma$ (see Figure~\ref{fig:gamma}).
        Consider a grid function $\xg$ defined on the discrete boundary $\gamma$. We can then define the \textit{difference potential with density $\xg$} as
        \begin{equation} \label{eq:Potential_def}
            \PN \xg \defeq w - \G^{(h)} \left( \Lh w {\big \rvert}_{\M^+} \right),\quad w = \begin{cases} \xg &\text{on } \gamma \\ 0 &\text{on } \N_0 \backslash \gamma \end{cases}
        \end{equation}
        The operation $\Lh w {\big \rvert}_{\M^+}$ in (\ref{eq:Potential_def}) represents first applying the operator $\Lh$ to the grid function $w$, then truncating the result to the grid set $\M^+$. The difference potential $\PN \xg$ is a grid function defined on $\N^+$ (hence the notation). It satisfies the homogeneous finite difference equation $\L^{(h)}(\PN \xg) = 0$ on $\M^+$. By truncating the difference potential to the grid boundary, we obtain the projection operator $P_{\gamma}$:
        \begin{equation} \label{eq:def_Pgamma}
            P_{\gamma} \xg \defeq (\PN \xg) {\big \rvert}_{\gamma} .
        \end{equation}
        The projection $P_{\gamma}$ defined by (\ref{eq:def_Pgamma}) has the following property: a grid function $\xg$ satisfies the difference Boundary Equation with Projection (BEP)
        \begin{equation} \label{eq:discrete_BEP}
            P_{\gamma} \xg + Tr^{(h)} \Gh g = \xg
        \end{equation}
        if and only if there is a solution $u$ on $\N^+$ to the finite difference equation \eqref{eq:aux_FD_scheme} such that $\xg$ is the trace of $u$ on the grid boundary $\gamma$. In this case,  $u$ is reconstructed by means of the discrete generalized Green's formula
        \begin{equation} \label{eq:potential_solution}
            u = \PN \xg + \Gh g
        \end{equation}
        In particular, the discrete right-hand side $g$ in equations (\ref{eq:discrete_BEP}) and (\ref{eq:potential_solution}) may be given by (\ref{eq:rhs_g}): $g=\Bh f$.
        Then,  the discrete BEP (\ref{eq:discrete_BEP}) equivalently reduces the fourth order accurate discrete approximation of the Helmholtz equation $\L u = f$ from the grid domain $\N^+$ to the grid boundary $\gamma$. 
        It will be convenient to specifically study the case where the governing equation is homogeneous, i.e. $f \equiv 0$. In this case, \eqref{eq:discrete_BEP} reduces to
        \begin{equation} \label{eq:discrete_BEP_homo}
            P_{\gamma} \xg = \xg
        \end{equation}
        Similar to \eqref{eq:discrete_BEP} and \eqref{eq:potential_solution}, solutions of \eqref{eq:discrete_BEP_homo} can be used to reconstruct the corresponding solution $u$ with the use of the difference potential
        \begin{equation} \label{eq:potential_solution_homo}
            u = \PN \xg
        \end{equation}

%--------------------------------------------%
%
%       Subsubsection:
%       Equation-Based Extension
%
%--------------------------------------------%

        \subsubsection{Equation-Based Extension} \label{sec:eq_based_ext}
        In order for $u$ from \eqref{eq:potential_solution} to approximate the solution of \eqref{eq:helmholtz} on $\N^+$, the grid density $\xg$ must be related, in a certain way, to the trace of the solution $u$ at the continuous boundary $\Gamma$. This relation is expressed by the \textit{extension operator}. Consider a pair of functions defined on $\Gamma$: $\xG = \left( \xi_0,\,\xi_1 \right) {\big\rvert}_\Gamma$. One can consider $\xi_0$ and $\xi_1$ as the Dirichlet and Neumann data, respectively, of some function $v = v(x,y)$ on $\Omega_0$:
        \begin{equation*}
            (\xi_0,\, \xi_1) {\big \rvert}_{\Gamma} = \left. \left( v,\, \frac{\partial v}{\partial {\bf n}} \right) \right|_{\Gamma}
        \end{equation*}
        This function $v$ can be defined in the vicinity of $\Gamma$ as a truncated Taylor expansion, with $\rho$ representing the distance (with sign) from the point of evaluation to $\Gamma$:
        \begin{equation} \label{eq:general_extension}
            v(x,y) \stackrel{\text{def}}{=} \left. v \right|_\Gamma + \rho \left. \frac{\partial v}{\partial {\bf n}} \right|_\Gamma + \frac{\rho^2}{2} \left. \frac{\partial^2 v}{\partial {\bf n}^2} \right|_\Gamma + \frac{\rho^3}{6} \left. \frac{\partial^3 v}{\partial {\bf n}^3} \right|_\Gamma + \frac{\rho^4}{24} \left. \frac{\partial^4 v}{\partial {\bf n}^4} \right|_\Gamma
        \end{equation}
        The definition (\ref{eq:general_extension}) of the new function $v(x,y)$ is not complete until the higher order normal derivatives are provided. These can be obtained using equation-based differentiation applied to the Helmholtz equation \eqref{eq:helmholtz}, where we assume $v$ is a solution and $v$ and $\partialdiff{v}{\bf n}$ are known analytically on $\Gamma$. When the domain $\Omega$ is a square, the outward normal derivatives on $\Gamma$ can be interpreted as standard $x-$ or $y-$ derivatives (or their negative counterparts), depending on which portion of the boundary one is considering.

        For example, let the right side of the square be $x=X=\text{const}$. Then, the outward normal derivative becomes the positive $x-$ derivative, by rearranging \eqref{eq:helmholtz}, we immediately get an expression for the second $x-$derivative evaluated along $\Gamma$:
        \begin{equation} \label{eq:second_normal}
            \frac{\partial^2 v}{\partial x^2}(X,y) = f(X,y) - \frac{\partial^2 v}{\partial y^2}(X,y) - k^2 v(X,y)
        \end{equation}
        In this arrangement, $v(X,y)$ can be replaced with the known $\xi_0(y)$, and $\frac{\partial^2 v}{\partial y^2}(X,y)$ can be replaced with its second tangential derivative, $\xi_0^{\prime \prime} (y)$. The third and fourth derivatives can also be obtained by first differentiating \eqref{eq:helmholtz} with respect to $x$, then subsequently replacing $v(X,y)$ with $\xi_0 (y)$,  $\partialdiff{v}{x} (X,y)$ with $\xi_1(y)$, and $\frac{\partial^2 v}{\partial x^2} (X,y)$ with the right-hand side of \eqref{eq:second_normal}. This process yields the following expressions:
        \begin{subequations} \label{eq:cartesian_extension}
            \begin{align}
                v(X,y) &= \xi_0(y) \label{eq:cart_ext_0}\\
                \frac{\partial v}{\partial x}(X,y) &= \xi_1(y)\label{eq:cart_ext_1}\\
                \frac{\partial^2 v}{\partial x^2}(X,y) &= f(X,y) - \xi_0^{\prime \prime} (y) - k^2 \xi_0(y) \label{eq:cart_ext_2}\\
                \frac{\partial^3 v}{\partial x^3}(X,y) &= \frac{\partial f}{\partial x} (X,y) - \xi_1^{\prime \prime} (y) - k^2 \xi_1(y) \label{eq:cart_ext_3}\\
                \frac{\partial^4 v}{\partial x^4}(X,y) &= \frac{\partial^2 f}{\partial x^2} (X,y) - \frac{\partial^2 f}{\partial y^2} (X,y) - k^2 f(X,y) + \xi_0^{(4)}(y) + 2k^2 \xi_0^{(2)} (y) + k^4 \xi_0(y) \label{eq:cart_ext_4}
            \end{align}
        \end{subequations}
        The expressions in \eqref{eq:cartesian_extension} can be substituted into \eqref{eq:general_extension} to calculate the values of $v(x,y)$ near the right side of $\Gamma$. Similar derivations can be used to compute $v(x,y)$ near other sides of the square, keeping in mind that the outward normal derivative on the left and bottom sides of the square correspond to the negative $x-$ and $y-$ derivatives, respectively.

        The function $v=v(x,y)$  can be constructed starting from any pair of functions $\left( \xi_0, \xi_1 \right)$ defined on $\Gamma$ by means of substituting \eqref{eq:cart_ext_0}-\eqref{eq:cart_ext_4} into the Taylor expansion \eqref{eq:general_extension}. Then, sampling $v$ only on the grid boundary $\gamma$, we define the extension operator ${\bf Ex}$ that yields the grid function $\xg$:
        \begin{equation*} %\label{eq:extension_def}
            \xg = {\bf Ex} (\xi_0, \xi_1) = v {\big \rvert}_\gamma .
        \end{equation*}
        As seen in \eqref{eq:cartesian_extension}, the operator $\Ex$ depends on the source term $f$. Hence, $\Ex$ is an affine operator:
        \begin{equation} \label{eq:extension_split}
            \Ex \xG = \Ex^{(H)} (\xi_0,\xi_1) + \Ex^{(I)} f
        \end{equation}
        where $\Ex^{(H)}$ represents its homogeneous (i.e., linear) part that only depends on $(\xi_0,\xi_1)$, and $\Ex^{(I)}$ is the inhomogeneous part that accounts for the source term from \eqref{eq:helmholtz}.

        Although the formulae for the normal derivatives \eqref{eq:cartesian_extension} were derived using the Helmholtz equation, $\xG = \left( \xi_0, \xi_1 \right)$ does not need to represent the Cauchy data of a solution $u$ to \eqref{eq:helmholtz} in order to apply the operator ${\bf Ex}$. However, if $\xG$ does correspond to a solution $u$: $\xG=\left. \left( u,\, \frac{\partial u}{\partial {\bf n}} \right) \right|_{\Gamma}$, then $\xg = \Ex \, \xG$ approximates  this solution near $\Gamma$ with fifth-order accuracy with respect to the grid size $h$, specifically at the grid nodes of $\gamma$.

        Let $u$ be a solution to \eqref{eq:helmholtz} on $\Omega$ in the homogeneous case, $f\equiv 0$, and let $\xG$ be the trace of $u$ along the continuous boundary $\Gamma$ such that $\xG = \left(u, \partialdiff{u}{\bf n} \right) {\big \rvert}_\Gamma$. Let $\xg = \Ex \xG$ and let $\PN \xg$ be the difference potential with density $\xg$. Let $p$ be the order of accuracy of the finite difference scheme (Sections \ref{sec:FDM} and \ref{sec:auxiliary_problem}). According to Reznik \cite{Reznik_83, reznik-english} (alternatively, see \cite{ryaben2002method}), as the grid $\N_0$ is refined, $\PN \xg$ converges to the solution $u$ (on the grid $\N^+$) with the convergence rate of $\mathcal{O}(h^p)$ provided that the number of terms in the Taylor expansion \eqref{eq:general_extension} is equal to $p+q$, where $q$ is the order of the differential operator $\L$. Given that the Helmholtz equation is second-order and we use a fourth-order finite difference scheme \eqref{eq:FD_scheme}, this would suggest the use of six terms in our extension. In practice, it has repeatedly been observed (see \cite{Medvinsky2012}, \cite{Medvinsky2015}, and \cite{Britt2013}, among others) that while sufficient, this bound is not tight, and the number of terms typically matches the order of accuracy of the finite difference scheme alone. Our use of four terms in \eqref{eq:general_extension} is corroborated by the numerical experiments in Section~\ref{sec:numerical}.

%--------------------------------------------%
%
%       Subsubsection:
%       Series Representation of the Boundary Data
%
%--------------------------------------------%

        \subsubsection{Series Representation of the Boundary Data} \label{sec:series_representation}
        Consider a set of basis functions, $\{\psi_j\}$, and the following two sets of pairs
        \begin{equation} \label{eq:basis_functions}
            {\bm \psi}^{(0)}_j = \left( \psi_j , 0 \right), \quad {\bm \psi}_j^{(1)} = \left( 0, \psi_j \right), \quad j = 1, ..., \infty
        \end{equation}
        Recall that we denote the boundary data  by $\xG = \left( \xi_0,\xi_1 \right)$, where $\xi_0$ represents the Dirichlet data and $\xi_1$ represents the Neumann data. Specifically, consider one smooth section of $\Gamma$ (i.e. one side of the square), denoted $\Gamma^*$, and let its boundary data be denoted $\xGs = \left( \xi_0^*, \xi_1^* \right)$. The separate components of this section of boundary data can be expanded individually along $\Gamma^*$:
        \begin{equation}\label{eq:boundary_inf}
            \xGs = \left( \xi_0^*, \xi_1^* \right) = \sum_{j=1}^{\infty} {\bf c}_j^{(0)} {\bm \psi}_{j}^{(0)} + \sum_{j=1}^{\infty} {\bf c}_j^{(1)} {\bm \psi}_{j}^{(1)} .
        \end{equation}
        The infinite series (\ref{eq:boundary_inf}) can be truncated after a finite number of terms to provide an approximation of $\xGs$. The number of terms $M^*$ is typically taken so as to make the truncated terms  negligible with respect to the accuracy attainable on the grid:
        \begin{equation} \label{eq:boundary_finite}
            \xGs = \sum_{j=1}^{M^*} {\bf c}_j^{(0)} {\bm \psi}_{j}^{(0)} + \sum_{j=1}^{M^*} {\bf c}_j^{(1)} {\bm \psi}_{j}^{(1)}
        \end{equation}
        Provided that the boundary data are sufficiently smooth, for the appropriately chosen  basis functions $\psi_j$ (e.g. Chebyshev, Fourier, etc...) the value of $M^*$ can be taken  relatively small.
        \begin{figure}[ht]
            \centering
            \includegraphics[width=0.2\textwidth]{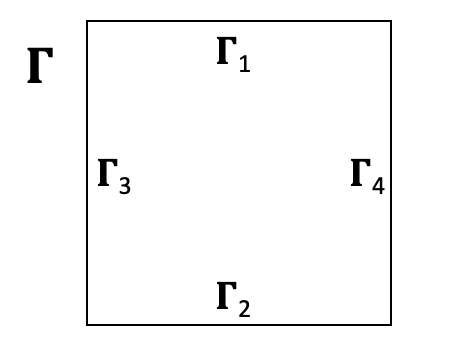}
            \caption{A labeling of the sides of $\Gamma$. The choice of ordering is arbitrary and only given here as a visual reference for the linear system in Section~\ref{sec:LinearSystem}.}
            \label{fig:gamma_decomp}
        \end{figure}
        The series representation  \eqref{eq:boundary_finite} can be extended to apply to all four sides of the square $\Gamma$ by combining the corresponding basis functions. Consider the labeling of the sides in Figure~\ref{fig:gamma_decomp}, and the following definition of the expanded set of basis functions $\Psi_j$:
        \begin{equation} \label{eq:Psi_def}
            \Psi_{j + (i-1)M^*} = \begin{cases} \psi_j &\text{on } \Gamma_i \\ 0 &\text{otherwise}  \end{cases} \quad \text{for $i = 1, ..., 4$} .
        \end{equation}
        Every element of $\Psi$ in (\ref{eq:Psi_def}) is defined on the entire $\Gamma$, while each $\xi_{\Gamma^i}$ has a series expansion independent of the others. Then, similar to \eqref{eq:basis_functions} we define the following pairs:
        \begin{equation*} %\label{eq:big_basis_functions}
            {\bm \Psi}^{(0)}_j = \left( \Psi_j , 0 \right), \quad {\bm \Psi}^{(1)}_j = \left( 0, \Psi_j \right), \quad j = 1, ..., M
        \end{equation*}
        where  $M = 4 \cdot M^*$, and  write the expansion of $\xG$ as
        \begin{equation} \label{eq:boundary_finite_gamma}
            \xG = \sum_{j=1}^{M} {\bf c}_j^{(0)} {\bm \Psi}_{j}^{(0)} + \sum_{j=1}^{M} {\bf c}_j^{(1)} {\bm \Psi}_{j}^{(1)}
        \end{equation}
        Note that, the choice of the same system of basis function for both the Dirichlet and Neumann data and for all four sides of the square is not a requirement, but it provides extra convenience for constructing the linear system in Section~\ref{sec:LinearSystem} and building the DDM extension in Section~\ref{sec:extension_2_domains}.

%--------------------------------------------%
%
%       Subsubsection:
%       Forming the Base Linear System
%
%--------------------------------------------%

        \subsubsection{Forming the Base Linear System} \label{sec:LinearSystem}
        Applying the extension operator (\ref{eq:extension_split}) to the series representation of $\xG$ in \eqref{eq:boundary_finite_gamma}, we have:
        \begin{align} \label{eq:extension_xiG}
            {\bf Ex}\, \xG &= \Ex^{(H)} \left( \sum_{j=1}^{M} {\bf c}_j^{(0)} {\bm \Psi}_{j}^{(0)} + \sum_{j=1}^{M} {\bf c}_j^{(1)} {\bm \Psi}_{j}^{(1)} \right) + \Ex^{(I)} f \nonumber \\
            &= \sum_{j=1}^{M} {\bf c}_j^{(0)} \Ex^{(H)} \, {\bm \Psi}_{j}^{(0)} + \sum_{j=1}^{M} {\bf c}_j^{(1)} \Ex^{(H)} \, {\bm \Psi}_{j}^{(1)} + \Ex^{(I)} f
        \end{align}
        Setting $\xg = \Ex \, \xG$ and substituting it into the BEP \eqref{eq:discrete_BEP} with $g=\Bh f$ yields:
        \begin{align*}
            &P_\gamma \xg = \xg - Tr^{(h)} \Gh \Bh f \\
            &P_\gamma \Ex \, \xG = \Ex \, \xG - Tr^{(h)} \Gh \Bh f \\
            &P_{\gamma} \left( \sum_{j=1}^{M} {\bf c}_j^{(0)} \Ex^{(H)} \, {\bm \Psi}_{j}^{(0)} + \sum_{j=1}^{M} {\bf c}_j^{(1)} \Ex^{(H)} \, {\bm \Psi}_{j}^{(1)} + \Ex^{(I)} f \right) \\
            & \hspace{3cm} = \sum_{j=1}^{M} {\bf c}_j^{(0)} \Ex^{(H)} \, {\bm \Psi}_{j}^{(0)} + \sum_{j=1}^{M} {\bf c}_j^{(1)} \Ex^{(H)} \bm \Psi_{j}^{(1)} + \Ex^{(I)} f - Tr^{(h)} \Gh \Bh f \\
            &\sum_{j=1}^{M} {\bf c}_j^{(0)} P_{\gamma} \Ex^{(H)} \, {\bm \Psi}_{j}^{(0)} + \sum_{j=1}^{M} {\bf c}_j^{(1)} P_{\gamma} \Ex^{(H)} \, {\bm \Psi}_{j}^{(1)}  + P_\gamma \Ex^{(I)} f \\
            & \hspace{3cm}= \sum_{j=1}^{M} {\bf c}_j^{(0)} \Ex^{(H)} \, {\bm \Psi}_{j}^{(0)} + \sum_{j=1}^{M} {\bf c}_j^{(1)} \Ex^{(H)} \, {\bm \Psi}_{j}^{(1)} + \Ex^{(I)} f - Tr^{(h)} \Gh \Bh f
        \end{align*}
        By collecting similar terms, we obtain the following system of linear algebraic equations
        \begin{align} \label{eq:linear_system_equations}
            &\sum_{j=1}^{M} {\bf c}_j^{(0)} (P_\gamma - I_\gamma) \Ex^{(H)} \, {\bm \Psi}_{j}^{(0)} + \sum_{j=1}^{M} {\bf c}_j^{(1)} (P_\gamma - I_\gamma) \Ex^{(H)} \, {\bm \Psi}_{j}^{(1)} \nonumber \\
            &\hspace{7cm}= (I_\gamma - P_\gamma) \Ex^{(I)} f - Tr^{(h)} \Gh \Bh f .
        \end{align}
        where $I_\gamma$ represents the identity operator in the space of grid functions $\xg$ defined on $\gamma$. System \eqref{eq:linear_system_equations} can be written in matrix form:
        \begin{equation} \label{eq:single_matrix_system}
            Qc = F
        \end{equation}
        where the matrix $Q = \left[Q^{(0)}, \, Q^{(1)} \right]$ is given by
        \begin{align} \label{eq:Q_construction}
            Q = {\Big [} & \underbrace{(P_\gamma - I_\gamma) {\bf Ex} \, {\bf \Psi}^{(0)}_1,\, ... \,,\, (P_\gamma - I_\gamma) {\bf Ex} \, {\bf \Psi}^{(0)}_M,}_{Q^{(0)}} \nonumber \\
            &\underbrace{(P_\gamma - I_\gamma) {\bf Ex} \, {\bf \Psi}^{(1)}_1,\, ... \,,\, (P_\gamma - I_\gamma) {\bf Ex} \, {\bf \Psi}^{(1)}_M}_{Q^{(1)}} {\Big ]} .
        \end{align}
        This matrix has dimension $|\gamma| \times 2M$ where $|\gamma|$ is the number of nodes in the grid boundary $\gamma$. The column vector
        \begin{equation}  \label{eq:c_construction}
            c = [\underbrace{c_1, \, ... \,,\, c_M}_{c^{(0) \intercal}}, \, \underbrace{c_{M+1},\, ... \,,\, c_{2M}}_{c^{(1) \intercal}}]^\intercal
        \end{equation}
        in equation (\ref{eq:single_matrix_system}) is a vector of unknowns with dimension $2M$, while  the vector $F$ has dimension $|\gamma|$ and  represents the inhomogeneous part of the problem: $(I_\gamma - P_\gamma) \Ex^{(I)} f - Tr^{(h)} \Gh \Bh f$. The first $M$ columns of $Q$ in (\ref{eq:Q_construction}) form the sub-matrix $Q^{(0)}$ and correspond to the coefficients $\{c^{(0)}_j\}$ in (\ref{eq:c_construction}), while columns $M+1$ through $2M$ form $Q^{(1)}$ and correspond to $\{c^{(1)}_j\}$.

        Note  that, solution to \eqref{eq:single_matrix_system} is not unique, as the system is derived only from  \eqref{eq:discrete_BEP}, (\ref{eq:rhs_g}) and does not take into account any boundary conditions. Therefore, we interpret the underdetermined system \eqref{eq:single_matrix_system}  as a core piece of the multi-subdomain decomposition algorithm, rather than a system to be solved in its own right. The decomposition algorithm is described in Section~\ref{sec:extension_2_domains} for the case of two subdomains and subsequently extended in Section~\ref{sec:extension_N_domains} to the case of a larger number of subdomains. For a discussion about implementing boundary conditions and completing the algorithm in the single domain case, see \cite{ryaben2002method, Medvinsky2012, Medvinsky2015, Britt_10, Britt2013}.

%--------------------------------------------%
%
%       Subsection:
%       Extension to 2 Subdomains
%
%--------------------------------------------%

        \subsection{Extension to 2 Subdomains} \label{sec:extension_2_domains}
        Reconsidering the problem of solving \eqref{eq:helmholtz_physical} over a partitioned domain as in Figure~\ref{fig:DDM_rect}. Let $\Gamma^{(i)}$ represent the boundary of $\Omega_i$, and let each $\Gamma^{(i)}$ be composed of its four sides as in Figure~\ref{fig:gamma_decomp}, so that $\Gamma^{(i,j)}$ denotes side $j$ of $\Gamma^{(i)}$. Further, define a new set of indices, $B$, to be the indices of $\Gamma^{(i,j)}$ that correspond to the boundary edges. For the two-domain case, this yields $B = \{(1,1), (1,2), (1,3), (2,1), (2,2), (2,4) \}$, as well as its complement $B^\complement = \{(1,4), (2,3)\}$ for the indices corresponding to both sides of the interface $\Sigma$. Let all grid sets and operators from Section~\ref{sec:base_subdomain} be defined for $\Omega_1$ and $\Omega_2$, independently. Partition and index the matrix $Q$ and the unknown column vector $c$ with the following notation:
        \begin{align*}
            Q^{(i,*,k)} &= \begin{bmatrix} Q^{(i,1,k)} &Q^{(i,2,k)} &Q^{(i,3,k)} &Q^{(i,4,k)} \end{bmatrix} \\
            c^{(i,*,k)} &= \begin{bmatrix} c^{(i,1,k) \intercal} &c^{(i,2,k) \intercal} &c^{(i,3,k) \intercal} &c^{(i,4,k) \intercal} \end{bmatrix}^\intercal .
        \end{align*}
        For $Q^{(i,j,k)}$, the indices $i\in \{1,2\}$ and $j \in \{1,2,3,4\}$ denote those columns corresponding to the basis functions defined to be non-zero over $\Gamma^{(i,j)}$. The index $k \in \{0,1\}$ distinguishes between the Dirichlet and Neumann data (compare to the notation $Q^{(0)}$ and $Q^{(1)}$ in Section~\ref{sec:LinearSystem}). The use of $c^{(i,j,k)}$ similarly identifies the coefficients of the corresponding basis functions over $\Gamma^{(i,j)}$ in either the Dirichlet or Neumann case. The independent linear systems for $\Omega_1$ and $\Omega_2$ can then be written as
        \begin{equation*}
            Q^{(1,*,*)} c^{(1,*,*)} = F^{(1)} \quad \text{and} \quad Q^{(2,*,*)} c^{(2,*,*)} = F^{(2)}
        \end{equation*}
        where $Q^{(i,*,*)} = \begin{bmatrix} Q^{(i,*,0)} &Q^{(i,*,1)} \end{bmatrix}$  and $c^{(i,*,*)} = \begin{bmatrix} c^{(i,*,0) \intercal} &c^{(i,*,1) \intercal} \end{bmatrix}^\intercal$. Note that the construction of $Q^{(i,*,*)}$ is identical to that of \eqref{eq:Q_construction} over a single subdomain. Equivalently, the independent linear systems can be expressed simultaneously as the block-diagonal system
        \begin{equation} \label{eq:2_block_system}
            \begin{bmatrix} Q^{(1,*,*)} &0 \\ 0 &Q^{(2,*,*)} \end{bmatrix} \begin{bmatrix} c^{(1,*,*)} \\ c^{(2,*,*)} \end{bmatrix} = \begin{bmatrix} F^{(1)} \\ F^{(2)} \end{bmatrix}
        \end{equation}
        Similar to \eqref{eq:single_matrix_system}, the solution to \eqref{eq:2_block_system} is not  unique because it is derived only from the discrete BEP \eqref{eq:discrete_BEP} combined with (\ref{eq:rhs_g}) for each subdomain and does not account for the boundary condition \eqref{eq:helmholtz_robin}. Additionally, since $\Omega$ has been decomposed into $\Omega_1$ and $\Omega_2$, an interface (or transmission) condition is needed to account for the lack of a boundary condition along the interface $\Sigma$.

%--------------------------------------------%
%
%       Subsubsection:
%       Boundary Conditions
%
%--------------------------------------------%

        \subsubsection{Boundary Conditions} \label{sec:boundary_conditions}
        To account for the boundary conditions, consider one $\Gamma^{(i,j)}$ (for $(i,j) \in B$). Substitute the series representation at the boundary (cf. \eqref{eq:boundary_finite}) into the boundary condition \eqref{eq:helmholtz_robin} for both $u$ and $\partialdiff{u}{\bf n}$. Expand the right-hand side of \eqref{eq:helmholtz_robin} as $\phi = \sum_m^{M^*} d^{(i,j)}_m \psi_m$ (using the same basis functions as in \eqref{eq:basis_functions}). Then,
        \begin{equation} \label{eq:BC_series}
            \alpha \left( \sum_{m=1}^{M^*} c^{(i,j,0)}_m \psi_m \right) + \beta \left( \sum_{m=1}^{M^*} c^{(i,j,1)}_m \psi_m \right) = \sum_m^{M^*} d^{(i,j)}_m \psi_m .
        \end{equation}
        Assuming that the basis functions $\psi_m$ are orthogonal, we derive from (\ref{eq:BC_series}):
        \begin{equation} \label{eq:BC_coefficients}
            \alpha c_m^{(i,j,0)} + \beta c_m^{(i,j,1)} = d_m^{(i,j)}, \quad \text{for } m \in \{1,...,M^*\}  .
        \end{equation}
        The $M^*$ equations (\ref{eq:BC_coefficients}) can be obtained for each index pair in $B$, adding a total of $6M^*$ extra equations. Note here that the sets of equations obtained for each $\Gamma^{(i,j)}$ are independent from one another, allowing greater flexibility in the boundary condition \eqref{eq:helmholtz_robin}. For example, the definitions of $\alpha$ and $\beta$ in \eqref{eq:helmholtz_robin} can be piece-wise constant, split along the different sections of $\Gamma$, i.e. on $\Gamma^{(i,j)}$
        \begin{equation}
            \alpha = \alpha^{(i,j)} \quad \beta = \beta^{(i,j)}
        \end{equation}
        where $\alpha^{(i,j)}$ and $\beta^{(i,j)}$ are constants and $\left( \alpha^{(i,j)}, \beta^{(i,j)} \right) \neq (0,0)$ for any pair $(i,j) \in B$. This generalization allows for both Dirichlet and Neumann conditions ($\beta^{(i,j)} = 0$ or $\alpha^{(i,j)} = 0$, respectively) as particular cases. The equations being added by \eqref{eq:BC_coefficients} are sparse compared to the rest of \eqref{eq:2_block_system}, which can be taken advantage of computationally, see Section~\ref{sec:solving_complete}. Further information on implementing mixed boundary conditions, as well as extending this process to include variable coefficient Robin conditions, can be found in \cite{Britt2013}.

%--------------------------------------------%
%
%       Subsubsection:
%       Transmission Condition
%
%--------------------------------------------%

        \subsubsection{Interface Condition} \label{sec:transmission_condition}
        The standard interface condition requires continuity of the solution and its flux  across the interface (see Section~\ref{sec:lions}). These two conditions can be enforced by equating the series representations of the Dirichlet data along $\Gamma^{(1,4)}$ and $\Gamma^{(2,3)}$, as well as setting the series representation for the Neumann data of $\Gamma^{(1,4)}$ equal to the negative of that for $\Gamma^{(2,3)}$.
        \begin{subequations} \label{eq:transmission_series}
        \begin{equation} \label{eq:dirichlet_transmission}
            \sum_{m=1}^{M^*} c_m^{(1,4,0)} \psi_m = \sum_{m=1}^{M^*} c_m^{(2,3,0)} \psi_m
        \end{equation}
        \begin{equation} \label{eq:neumann_transmission}
            \sum_{m=1}^{M^*} c_m^{(1,4,1)} \psi_m = - \sum_{m=1}^{M^*} c_m^{(2,3,1)} \psi_m .
        \end{equation}
        \end{subequations}
        As with the boundary conditions in Section~\ref{sec:boundary_conditions}, the use of identical sets of orthogonal basis functions along each side is exploited to obtain the following two sets of equations for each $m \in \{1,...,M^*\}$:
        \begin{subequations} \label{eq:transmission_coefficients}
        \begin{equation} \label{eq:transmission_coefficients_D}
            c_m^{(1,4,0)} - c_m^{(2,3,0)} = 0
        \end{equation}
        \begin{equation} \label{eq:transmission_coefficients_N}
            c_m^{(1,4,1)} + c_m^{(2,3,1)} = 0 .
        \end{equation}
        \end{subequations}
        Collectively, \eqref{eq:transmission_coefficients} provides $2M^*$ equations  to supplement \eqref{eq:2_block_system}.

        Alternative interface conditions can be chosen and implemented in a similar fashion. For example, if $u_i$ is the solution to the subproblem on $\Omega_i$, then for constants $a^{(0)}, a^{(1)}, b^{(0)}, b^{(1)}$ and smooth functions $\eta^{(0)}, \eta^{(1)}$, a class of interface conditions can be defined as follows on the interface $\Gamma^{(1,4)} = \Gamma^{(2,3)}$:
        \begin{equation} \label{eq:transmission_general}
           a^{(0)} u_1 + b^{(0)} u_2 = \eta^{(0)} \hspace{1.5cm} a^{(1)} \partialdiff{u_1}{\bf n_1} + b^{(1)} \partialdiff{u_2}{\bf n_2} = \eta^{(1)} .
        \end{equation}
        Any transmission conditions of type (\ref{eq:transmission_general}) can be accounted for by following the same steps as in \eqref{eq:transmission_series} and \eqref{eq:transmission_coefficients}. The only addition is to let $\eta^{(0)} = \sum_{m=1}^{M^*} \eta^{(0)}_m \psi_m$ and $\eta^{(1)} = \sum_{m=1}^{M^*} \eta^{(1)}_m \psi_m$ be the expansions of $\eta^{(0)}$ and $\eta^{(1)}$. This yields $a^{(0)} c_m^{(1,4,0)} + b^{(0)} c_m^{(2,3,0)} = \eta^{(0)}_m$ and $a^{(1)} c_m^{(1,4,1)} + b^{(1)} c_m^{(2,3,1)} = \eta^{(1)}_m$ as the conditions for the coefficients, where the choice of $a^{(0)} = a^{(1)} = b^{(1)} = 1$, $b^{(0)} = -1$, and $\eta^{(0)} \equiv \eta^{(1)} \equiv 0$ recovers the original condition. By allowing linear combinations and inhomogeneities in the interface conditions, a wider set of situations such as jumps over the interface in the solution, its flux, or both can be accounted for. In this paper, however, we are only considering the case where the solution and its flux are continuous on $\Sigma$.

%--------------------------------------------%
%
%       Subsubsection:
%       Solving the Complete System
%
%--------------------------------------------%

        \subsubsection{Solving the Complete System} \label{sec:solving_complete}
        By supplementing the system \eqref{eq:2_block_system} with the equations derived in \eqref{eq:BC_coefficients} and \eqref{eq:transmission_coefficients}, the complete system can be expressed as a new matrix equation
        \begin{equation} \label{eq:complete_Q}
            \overline{Q} c^{(*,*,*)} = \overline{F}
        \end{equation}
        where the dimension of $\overline{Q}$ is $(2|\gamma| + 8M^*) \times 16M^*$. The system \eqref{eq:complete_Q} can be solved by minimizing the $\ell_2$ norm $\|\overline{Q} c^{(*,*,*)} - \overline{F}\|_2$ through traditional least squares methods, e.g., a QR-factorization. Note that while the least squares solution is unique, it is the existence of a classical solution to \eqref{eq:helmholtz_physical} --- from which \eqref{eq:complete_Q} is ultimately derived --- that guarantees $\|\overline{Q} c^{(*,*,*)} - \overline{F}\|_2$ will be within discretization error of zero. In fact, if $M^*$ is chosen large enough in \eqref{eq:boundary_finite}, then $\|\overline{Q} c^{(*,*,*)} - \overline{F}\|_2$ decreases at a rate of $\mathcal{O}(h^4)$ (the order of accuracy of the finite difference scheme) as the grid is refined.

        Rather than adding equations \eqref{eq:BC_coefficients} and \eqref{eq:transmission_coefficients} to the system, these conditions can instead be resolved through substitution and the elimination of unknowns. For the boundary condition equations \eqref{eq:BC_coefficients}, first consider the case where \eqref{eq:helmholtz_robin} reduces to a Dirichlet boundary condition (i.e. $\alpha=1$, $\beta=0$). In this case, the coefficients $c^{(i,j,0)}$ (for $(i,j) \in B$) are obtained directly when expanding the right-hand side of \eqref{eq:helmholtz_robin}, eliminating those coefficients from the larger linear system. The coefficients in $c^{(i,j,0)}$ are multiplied by the corresponding columns of $Q^{(i,j,0)}$, then subtracted over to the right-hand side of \eqref{eq:2_block_system}. If \eqref{eq:helmholtz_robin} reduces to a Neumann condition (i.e. $\alpha=0$, $\beta=1$), the same process is followed but for $c^{(i,j,1)}$ and $Q^{(i,j,1)}$. In either case, this process eliminates $6M^*$ unknowns from the system ($M^*$ unknowns for each $(i,j) \in B$) leaving $10M^*$ unknowns rather than the original $16M^*$ unknowns.

        When \eqref{eq:helmholtz_robin} does not reduce to a Dirichlet or Neumann condition ($\alpha \neq 0$ and $\beta \neq 0$), we can still eliminate unknowns by means of substitution. Consider \eqref{eq:BC_coefficients}, and rearrange the terms to solve for either $c_m^{(i,j,0)}$ or $c_m^{(i,j,1)}$:
        \begin{equation} \label{eq:BCs_rewritten}
            c_m^{(i,j,1)} = \frac{1}{\beta} d_m^{(i,j)} - \frac{\alpha}{\beta} c_m^{(i,j,0)}
        \end{equation}
        From (\ref{eq:BCs_rewritten}), the $\frac{1}{\beta} d^{(i,j)}_m$ terms can be multiplied by the corresponding columns of $Q^{(i,j,1)}$ and subtracted to the right-hand side, while the $\frac{\alpha}{\beta} c^{(i,j,0)}_m$ terms can be combined with their like terms from the original system \eqref{eq:2_block_system}. Similar to the Dirichlet and Neumann cases, $6M^*$ unknowns are eliminated from the system.

        The interface conditions \eqref{eq:transmission_coefficients_D}
        can be accounted for by adding the respective columns, $Q^{(1,4,0)}_m$ and $Q^{(2,3,0)}_m$, and eliminating one of the coefficients, $c^{(1,4,0)}_m$ or $c^{(2,3,0)}_m$. As these conditions exist for $m \in \{1,...,M^*\}$, resolving the interface conditions this way eliminates $M^*$ unknowns from the system. Following the same process for \eqref{eq:transmission_coefficients_N} (subtracting columns instead of adding) eliminates an additional $M^*$ unknowns.

        In the case where \eqref{eq:BC_coefficients} and \eqref{eq:transmission_coefficients} are included as supplemental equations, the overall system has dimension ($2|\gamma| + 8M^*) \times 16M^*$. If the conditions are resolved, the dimension is $2|\gamma| \times 8M^*$, which enables faster solution.  The solution vector $c^{(*,*,*)}$ is used to reconstruct $\xi_{\Gamma^{(1)}}$ and $\xi_{\Gamma^{(2)}}$ through the series representation \eqref{eq:boundary_finite} for each subproblem. In turn, $\xi_{\Gamma^{(1)}}$ and $\xi_{\Gamma^{(2)}}$ are extended to their respective grid boundaries, as described in Section~\ref{sec:eq_based_ext}. Finally, a fourth-order accurate approximation to the unique solution of \eqref{eq:helmholtz_physical} is obtained by applying \eqref{eq:potential_solution} to the resulting $\xi_{\gamma^{(1)}}$ and $\xi_{\gamma^{(2)}}$. These approximations collectively provide an approximation of the global solution to \eqref{eq:helmholtz_physical} on the overall domain $\Omega$.

%--------------------------------------------%
%
%       Subsection:
%       Extension to N Subdomains
%
%--------------------------------------------%

        \subsection{Extension to N Subdomains} \label{sec:extension_N_domains}
\begin{figure}[ht]
            \centering
            \includegraphics[align=c,width=.3\textwidth]{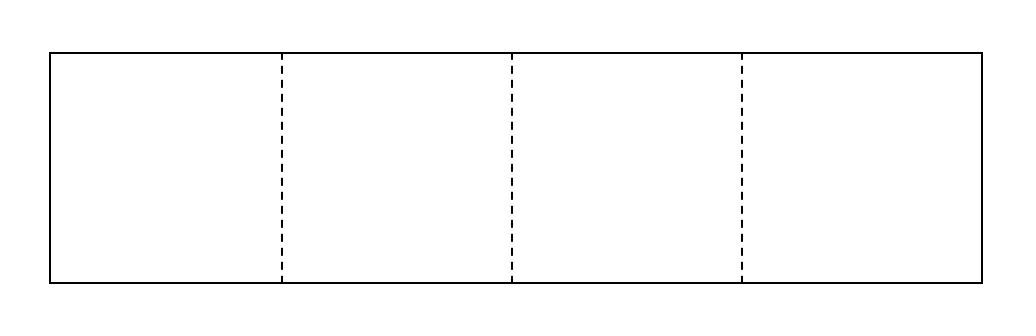}
            \includegraphics[align=c,width=.3\textwidth]{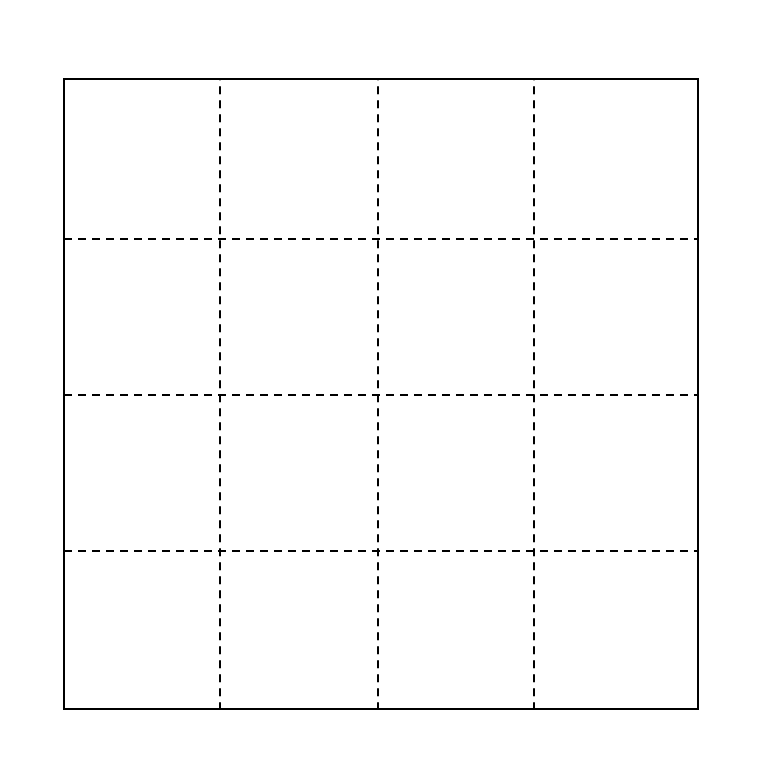}
            \includegraphics[align=c,width=.3\textwidth]{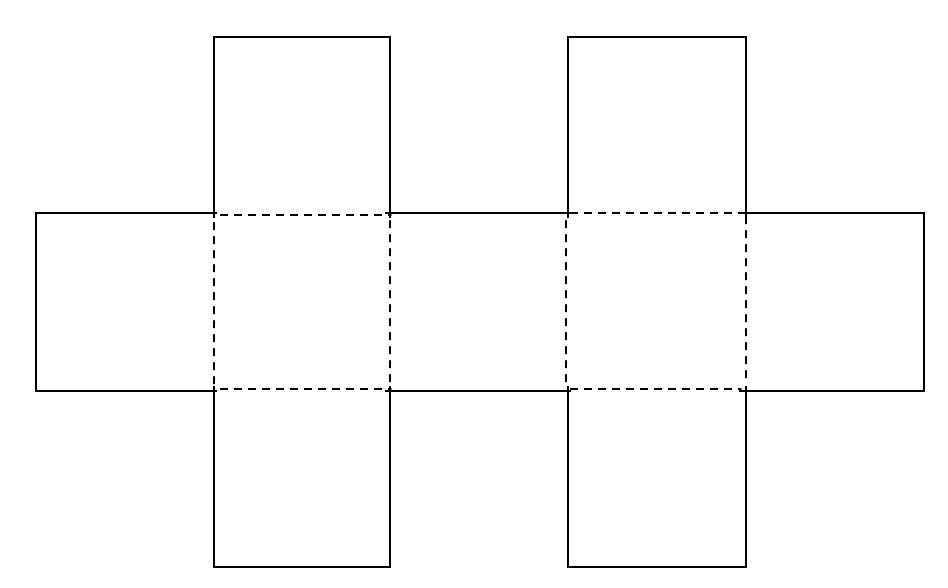}
            \caption{Examples of domains that have a valid $N-$subdomain decomposition.}
            \label{fig:N_domains_valid}
        \end{figure}

        The extension to $N$ subdomains is a natural extension of the two-subdomain case. Consider \eqref{eq:helmholtz_physical} over a domain $\Omega$ that is split into $N$ identical (square) subdomains, whose interfaces are full edges of the squares (see Figure~\ref{fig:N_domains_valid}).
        Returning to the triple index notation used in Section~\ref{sec:extension_2_domains}, let the first argument vary from $1$ to $N$, rather than stopping at $2$, and let all grid sets and operators from Section~\ref{sec:base_subdomain} be defined independently for each $\Omega_i$. To build the matrix for the linear system, combine the $Q^{(i,*,*)}$ from each subdomain in a block-diagonal style. The  vectors of unknowns and  right-hand sides  from each subdomain are simply appended to create the following system:
        \begin{equation} \label{eq:N_block_system}
            \begin{bmatrix} Q^{(1,*,*)} &0 &0 &0 \\ 0 &Q^{(2,*,*)} &0 &0 \\ 0 &0 &\ddots &0 \\ 0 &0 &0 &Q^{(N,*,*)} \end{bmatrix} \begin{bmatrix} c^{(1,*,*)} \\ c^{(2,*,*)} \\ \vdots \\ c^{(N,*,*)} \end{bmatrix} = \begin{bmatrix} F^{(1)} \\ F^{(2)} \\ \vdots \\ F^{(N)} \end{bmatrix}
        \end{equation}
        To generalize the handling of boundary and interface conditions, extend the definition of the set $B$
        \begin{equation}
            B = \left\{ (i,j) {\big \vert} \Gamma^{(i,j)} \cap \dO \neq \emptyset \right\}
        \end{equation}
        so that $|B \cup B^\complement| = 4N$. If $(i,j) \in B$, then $\Gamma^{(i,j)}$ has an associated boundary condition specified by \eqref{eq:helmholtz_robin} and the process described in Section~\ref{sec:boundary_conditions} can be applied for each $(i,j) \in B$ to obtain the necessary supplemental equations. If $(i,j) \in B^\complement$, then $\Gamma^{(i,j)}$ is an interface, requiring the process from Section~\ref{sec:transmission_condition} to determine the supplemental equations. Adding these equations yields the $N$ subdomain version of \eqref{eq:complete_Q}:
        \begin{equation} \label{eq:N_domain_complete}
            \overline{Q}_{N} c^{(*,*,*)} = \overline{F}_{N} .
        \end{equation}
        where $\overline{Q}_{N}$ and $\overline{F}_{N}$ represent the matrix from the left-hand side of \eqref{eq:N_block_system} and the vector of the right-hand side, respectively, after being supplemented with boundary and interface condition equations. The shape of the domain determines how many equations correspond to boundary conditions as opposed to interface conditions, but there will always be $4N M^*$ equations added to \eqref{eq:N_block_system} ($M^*$ equations for each $\Gamma^{(i,j)}$). As in Section~\ref{sec:solving_complete}, these equations can often be resolved with substitution and elimination to reduce the cost of solving the linear system.

%--------------------------------------------%
%
%       Subsection:
%       Implementation Details
%
%--------------------------------------------%

    \subsection{Implementation Details} \label{sec:implementation}
    In this section, we provide the important implementation details of the proposed algorithm, which are further justified in Section~\ref{sec:complexity}.  We assume that all subdomains are identical squares and that $\Omega$ satisfies the requirements described at the beginning of Section~\ref{sec:extension_N_domains} and in Figure~\ref{fig:N_domains_valid}. In \eqref{eq:helmholtz}, we assume that the wavenumber $k$ is  piece-wise constant over $\Omega$, and constant on any given $\Omega_i$. Such assumptions allow the exploration of the best-case scenario. A brief discussion of possible generalizations is given in Section~\ref{sec:conclusions}. Consider the following summary of the algorithm:
    \begin{enumerate}
        \item For each $\Omega_i$:
        \begin{enumerate}
            \item Define the auxiliary problem \eqref{eq:auxiliary_problem}, as well as the grid sets and operators from Section~\ref{sec:grid_sets_and_operators}.
            \item Compute $Q^{(i,*,*)}$, the left-hand side of \eqref{eq:single_matrix_system}.
            \item Compute $P_\gamma {\bf Ex}^{(I)} f$ and $Tr^{(h)} \Gh \Bh f$ to form $F^{(i)}$ from the right-hand side of \eqref{eq:linear_system_equations}.
        \end{enumerate}
        \item Solve \eqref{eq:N_domain_complete} for $c^{(*,*,*)}$:
        \begin{enumerate}
            \item Assemble \eqref{eq:N_block_system} and resolve the boundary conditions and interface conditions (either as supplemental equations or as in Section~\ref{sec:solving_complete}).
            \item Compute the QR-factorization of $\overline{Q}_{N} = QR$.
            \item Compute $c^{(*,*,*)} = R^{-1} Q^* \overline{F}_{N}$ where $Q^*$ is the conjugate transpose of $Q$.
        \end{enumerate}
        \item For each $\Omega_i$:
        \begin{enumerate}
            \item Use $c^{(i,*,*)}$ and the series representation \eqref{eq:boundary_finite} for each $\Gamma^{(i,j)}$ to reconstruct the boundary data.
            \item Extend $\xi_{\Gamma^{(i)}}$ using ${\bf Ex}$ to obtain $\xi_{\gamma}$.
            \item Apply \eqref{eq:potential_solution} to $\xi_{\gamma}$ to obtain a local, fourth-order approximation to the solution $u$ on $\Omega_i$.
        \end{enumerate}
    \end{enumerate}
    The local solutions are then assembled to collectively provide a global approximation of the solution $u$ on $\Omega$. Note, that the entirety of steps 1 and 3 can be distributed on parallel processors for each subdomain. Once the algorithm has been run, the structure of the method allows several problem variations to be solved more economically because they do not affect the structure of terms that have already been computed in specific parts of steps 1 and 2.

    The first part of the algorithm that requires special consideration is step 1(b). Constructing $Q^{(i,*,*)}$ is expensive (see Section~\ref{sec:complexity}), but in general we do not need to recompute it every time the algorithm is run. Due to the use of geometrically identical subdomains as well the same set of basis functions throughout the problem, the only factor that distinguishes $Q^{(i,*,*)}$ from $Q^{(j,*,*)}$ is the wavenumber $k$ on $\Omega_i$ and $\Omega_j$. Therefore, $Q^{(i,*,*)}$ can be reused for any subdomain $\Omega_j$ such that the value of $k$ is shared across both $\Omega_i$ and $\Omega_j$. The cost to construct $Q^{(i,*,*)}$ should only be accrued once for each unique value of $k$ across all subdomains. In the case where $k$ is uniform across $\Omega$, $Q^{(i,*,*)} = Q^{(j,*,*)}$ for all $i,j \in \{ 1,...,N \}$, so the base linear system is only computed once, regardless of the number of subdomains. Further, as long as each $Q^{(i,*,*)}$ is saved after being computed, it can be reused in future problems for subdomains with the corresponding value of $k$, thus allowing the algorithm to run without constructing any $Q^{(i,*,*)}$ matrices. In this sense, we consider $Q^{(i,*,*)}$ to be pre-computed, thereby separating the cost of its construction from the run-time complexity of the algorithm.

    In step 2(b), a QR factorization is used to find the least squares solution of the matrix equation \eqref{eq:N_domain_complete}. The cost of QR factorization grows as the number of subdomains increases (see Section~\ref{sec:complexity}). However, once the factorization has been performed, changes to the right-hand sides of \eqref{eq:helmholtz} and \eqref{eq:helmholtz_robin} (i.e., $f$ and $\phi$, respectively) do not affect the left-hand side of \eqref{eq:N_domain_complete}. Thus, for a series of problems where only $f$ and $\phi$ vary, the cost of the QR factorization is only accrued on the first problem, effectively sharing its cost between such problems. Examples of the time saved in such cases are reported in Section~\ref{sec:numerical}. Further, in the case where $\phi$ changes while $f$ remains the same, step 1(c) can also be reused, thus starting the algorithm from step 2(c) and saving the cost of applying $\Gh$ in step 1(c).

    Finally, if the type of boundary condition is changed on a given $\Gamma^{(i,j)}$ by changing the piecewise-constant values of $\alpha$ or $\beta$ on the left-hand side of \eqref{eq:helmholtz_robin}, then the algorithm can begin at step 2(a). However, in practice we do not exploit this case for time savings as step 2(a) is relatively inexpensive to compute.

%--------------------------------------------%
%
%       Subsection:
%       Complexity
%
%--------------------------------------------%

    \subsection{Complexity} \label{sec:complexity}
    The complexity of the algorithm depends on two main factors: Solving the discrete AP (i.e. applying the operator $\Gh$) and computing the QR-factorization of $\overline{Q}_N$ from \eqref{eq:N_domain_complete}. Further, the applications of $\Gh$ include the pre-computed construction of $Q^{(i,*,*)}$ and the run-time steps 1(c) and 3(c) of the algorithm. We emphasize that our algorithm provides the exact solution of the discrete Helmholtz problem, as opposed to the traditional DDMs that are typically iterative. Due to the non-iterative nature of our method, a direct comparison of its complexity to that of  the conventional DDMs is poorly defined and  not explored in this paper. We, however, provide a thorough analysis of the complexity of our method as it depends on the various parameters of the discretization.

    First, consider applications of $\Gh$. This operator is only applied to individual subdomains, so let $n$ be the number of grid nodes in one direction in the discretization of an auxiliary domain. Recall from Section~\ref{sec:auxiliary_problem} the choice of boundary conditions for the $y-$boundaries in \eqref{eq:auxiliary_problem} and the requirement that $k$ be constant on $\Omega_i$. These choices allow the discrete AP to be solved with a combination of a sine-FFT in the $y-$direction and a tridiagonal solver in the $x-$direction, yielding a complexity of $\mathcal{O}(n^2 \log n)$. This is the contribution of one application of $\Gh$, but $\Gh$ gets applied many times over the course of the method. In particular, the construction of any $Q^{(i,*,*)}$ requires $8M^*$ applications of $\Gh$ --- one for each column --- giving the construction of  $Q^{(i,*,*)}$ a complexity of $\mathcal{O} \left( M^* n^2 \log n \right)$. In the worst-case scenario where every subdomain has a unique value of $k$, the construction of all $N$ of the distinct $Q^{(i,*,*)}$ matrices is $\mathcal{O} \left( NM^* n^2 \log n \right)$. However, it is important to note that the columns of $Q^{(i,*,*)}$ are  independent of one another, allowing the construction of $Q^{(i,*,*)}$ to be distributed to a number of parallel processors up to the horizontal dimension of the matrix. Furthermore, $\Gh$ is also applied twice to every subdomain to construct the right-hand side of \eqref{eq:N_block_system}, and an additional application is required in \eqref{eq:potential_solution} to obtain the final approximation. These $3N$ applications contribute $\mathcal{O} (N n^2 \log n)$ to the overall method complexity.

    The cost of the QR-factorization of $\overline{Q}_N$ from \eqref{eq:N_domain_complete} is the other main contribution to the method's complexity. Assuming that \eqref{eq:N_domain_complete} is formed by resolving the boundary and transmission conditions as discussed in Section~\ref{sec:solving_complete}, the dimensions of $\overline{Q}_N$ are $N |\gamma| \times 4NM^*$. Note that $|\gamma|$ is roughly proportional to $n$ because $\gamma$ only contains those nodes closest to the boundary of the subdomain, see Figure~\ref{fig:gamma}. The QR-factorization of a matrix depends linearly on the first dimension and quadratically on the second, giving our QR-factorization a complexity of $\mathcal{O} \left( (Nn) (NM^*)^2 \right)$, or equivalently, $\mathcal{O} \left(N^3 n (M^*)^2 \right)$. Further, $M^*$ is typically held constant for a given collection of problems (the selection of $M^*$ is discussed further in Section~\ref{sec:numerical}), so we consider the cost to be $\mathcal{O} \left(N^3 n \right)$. When the number of subdomains is large, it will become necessary to avoid repeating the QR-factorization when possible as discussed in Section~\ref{sec:implementation}, which leaves us with the $\mathcal{O} \left( N^2 n \right)$ operation of multiplying $R^{-1}Q^*$ by the source vector $\overline{F}_N$.

    Combining the costs of solving the AP and computing the QR-factorization gives the method an overall complexity of $\mathcal{O} \left( N^3 n + N n^2 \log n \right)$. For comparison, consider a simplified situation: let $\Omega$ be a square, and let $N = N_d^2$ such that there are $N_d$ subdomains along each side of $\Omega$, allowing the complexity to be rewritten as $\mathcal{O} \left( n N_d^6 + n^2 N_d^2 \log(n) \right)$. Consider solving \eqref{eq:helmholtz_physical} over this $\Omega$ with a finite-difference method and without domain decomposition. If the wavenumber is constant and uniform, and the boundary conditions on either the $x-$ or $y-$boundaries are homogeneous Dirichlet conditions, then we can directly use the FFT-based solver mentioned earlier. This domain has $nN_d$ nodes in each direction, so the complexity of this method would be $\mathcal{O} \left( (n N_d)^2 \log(n N_d) \right)$. This complexity is better than that of the proposed  method. Yet we stress that the FFT-based solver is only applicable in this simplest case, and is inflexible in terms of boundary conditions and variation of the wavenumber. In order to relax the requirements on the wavenumber and boundary conditions, we would need to resort to an LU or similar factorization to invert the finite difference operator with a complexity of $\mathcal{O} \left( ((nN_d)^2)^3 \right) = \mathcal{O} \left( (nN_d)^6 \right) $. This approach can capture a wide range of boundary conditions and a variable wavenumber, but it is still ill-suited to cases with piecewise-constant $k$, as the global solution loses regularity at the interfaces. In these simple cases where all methods apply, the FFT- and LU-based solvers provide lower and upper bounds for the expected performance of our method. However, unlike these two methods, our method extends naturally to more complex domain shapes and boundary conditions. 

%--------------------------------------------%
%
%       Section:
%       Numerical Results
%
%--------------------------------------------%

\section{Numerical Results} \label{sec:numerical}
In this section, we present numerical results corroborating the fourth-order convergence of the method, as well as the theoretical computational costs as discussed in Section~\ref{sec:complexity}. For all the test cases, we consider the Helmholtz equation \eqref{eq:helmholtz_physical} where the wavenumber $k$ is constant on any given $\Omega_i$ but piecewise constant over $\Omega$. The coefficients $\alpha$ and $\beta$ from \eqref{eq:helmholtz_robin} are piecewise constant on each $\Gamma^{(i)}$ such that each is constant on any given subdomain edge, $\Gamma^{(i,j)}$. However, for convenience in presentation, most examples use a uniform boundary condition across the entire $\dO$.

In all the test cases, we choose the domain $\Omega$ such that it can be split into  $N$ identical, square subdomains $\Omega_i$, whose interfaces are all full edges. Each $\Omega_i$ has side length $2$ and every corresponding auxiliary domain is a square with side length $2.2$. For consistency, we always let $\Omega_1$ be centered at the origin. Two particular domain shapes that provide convenient and systematic settings for analysis are a long duct and a large square. The duct is a quasi-one dimensional decomposition where $\Omega_{i+1}$ extends from $\Omega_i$ in the positive $x-$direction, allowing us to directly observe various behaviors of the method with respect to the number of subdomains $N$. On the other hand, a square domain will be decomposed into $N = N_d^2$ subdomains as discussed in Section~\ref{sec:complexity}, where $\Omega_1$ is again centered at the origin, and acts as the ``bottom-left" corner of $\Omega$ (see Figure~\ref{fig:checkerboard_square}), with $N-1$ subdomains extending in both the positive $x-$ and $y-$ directions. This domain gives us less direct control over $N$ itself, but it provides a framework to observe the method's performance in the presence of an increasing number of interior cross-points.

All auxiliary problems are solved with the method of difference potentials employing the fourth-order accurate compact finite difference scheme \eqref{eq:FD_scheme} on a series of Cartesian grids, starting with $n = 64$ cells uniformly spaced in each direction and progressively doubling $n$ with each refinement. The number of basis functions $M^*$ used in the expansion of $\xi_{\Gamma^{(i,j)}}$ is generally chosen grid-independent \cite{Medvinsky2012}, such that the boundary data are represented to a specified tolerance that is smaller than the error attainable on the given grids. Further increasing $M^*$ offers little to no benefit in the final accuracy of the method as we are still limited by the accuracy of the finite difference scheme. It has even been observed that selecting $M^*$ too large can have adverse effects on the overall accuracy \cite{Britt2013}, particularly on coarse grids. In this event, it is sufficient to simply reduce $M^*$ for the coarse grids, a practice that we indicate in the relevant results.

There are two kinds of test problems in this section: those with a known exact solution, and those without a known solution. For the test cases with a known solution, the source term and boundary data are derived by substituting the solution into the left-hand side of \eqref{eq:helmholtz} and \eqref{eq:helmholtz_robin}, respectively. In this case, the error is computed by taking the maximum norm of the difference between the approximated and the exact solution on the grid $\M^+$ (across all subdomains). The convergence rate is then determined by taking the binary logarithm of the ratio of errors on successively refined grids. In general, these test cases either have a uniform wavenumber throughout the domain, or are posed across a small number of subdomains in order to simplify the derivation of an exact solution.

On the other hand, when we want to specify the boundary conditions, source function, or piecewise constant wavenumber, we do not necessarily have an exact solution available and therefore  cannot compute the error directly to determine convergence. Instead, we introduce a grid-based metric where we compare the approximations on the shared nodes of successively refined grids. For a grid with $n\times n$ nodes, we denote the corresponding approximation by $\uhn$. Because of how the grids are structured, the nodes of the $\frac{n}{2} \times \frac{n}{2}$ grid are a subset of those in the $n \times n$ grid, so we can compute the maximum norm of the difference between these approximations, $\|\uhn - \uhnh\|_{\infty}$, on the nodes of $\M^+$ from the $\frac{n}{2} \times \frac{n}{2}$ grid. Similar to the first case, we can then estimate the convergence rate by considering the binary logarithm of the ratio of these maximum norm differences on successive grids. The new convergence metric does not evaluate the actual error. If, however, the discrete solution converges to the continuous one with a certain rate in the proper sense, then this alternative metric will also indicate convergence with the same rate regardless of whether the continuous solution is known or not and so it is convenient to use when the true solution is not available.

All of the computations in this section were performed in MATLAB (ver. R2019a) and used the package Chebfun \cite{chebfun} to handle all Chebyshev polynomial related operations. The QR-factorizations were performed using MATLAB's built-in `economy-size' QR-factorization.

%--------------------------------------------%
%
%       Subsection:
%       Uniform Wavenumber
%
%--------------------------------------------%

    \subsection{Uniform Wavenumber} \label{sec:uniform_k}
    In order to measure the complexity of the solver, we start with the case of one subdomain (i.e. $\Omega_1 \cong \Omega$). We use the homogeneous test solution $u = e^{i\frac{k}{\sqrt{2}}{(x+y)}}$ with wavenumber $k=13$, $M^*=40$, and Robin boundary conditions defined by $\alpha = 1$ and $\beta=1$ in \eqref{eq:helmholtz_robin}. The results in Table \ref{tab:base_case} corroborate both the fourth-order convergence rate of the overall method and the computational complexity of solving the discrete AP. As discussed in Section~\ref{sec:complexity}, for a grid with $n$ nodes in each direction, the FFT-based solver should have a complexity of $\mathcal{O}(n^2 \log n ) $, which produces the scale factors of approximately $4$ in the right-most column of Table \ref{tab:base_case} as $n$ is doubled. Note that Table \ref{tab:base_case} also corroborates the same complexity for the construction of $Q$ from \eqref{eq:single_matrix_system} (or equivalently, $Q^{(i,*,*)}$ from \eqref{eq:N_block_system}) because the dominating cost of constructing $Q$ is the application of $\Gh$ for every basis function in the subdomain. Note  that, the $\Gh$ timings will remain approximately constant for any number of subdomains $N$ (up to the available number of processors), with only small increases due to the overhead incurred by parallel communication.

    \begin{table}[ht]
        \centering
        \begin{tabular}{r c l l c l l}
            \toprule
            $n$ & &Error &Rate & &$\Gh$ time &Ratio \\ \cmidrule{3-4} \cmidrule{6-7}
            64 & &1.05e-03 &- & &0.0064 &- \\
            128 & &6.42e-05 &4.04 & &0.0083 &1.30 \\
            256 & &3.94e-06 &4.03 & &0.042 &5.07 \\
            512 & &2.47e-07 &4.00 & &0.186 &4.44 \\
            1024 & &1.53e-08 &4.01 & &0.824 &4.44 \\
            2048 & &9.89e-10 &3.95 & &3.430 &4.17 \\
            \bottomrule
        \end{tabular}
        \caption{Grid convergence and time (in seconds) to apply $\Gh$ to solve the discrete AP for the single subdomain base case. The test solution is $e^{i\frac{k}{\sqrt{2}}{(x+y)}}$ with $k=13$, the Robin boundary conditions are defined by $\alpha=\beta=1$ in \eqref{eq:helmholtz_robin}, and $M^*=40$.}
        \label{tab:base_case}
    \end{table}

    Next, we consider cases where $\Omega$ can be decomposed into two or more subdomains, under the simplifying assumption that $k$ is uniformly constant throughout $\Omega$, sharing the same value on every $\Omega_i$. This case is the simplest to consider because every $\Omega_i$ will use the same $Q^{(i,*,*)}$, removing the need to compute a new $Q^{(i,*,*)}$ for every possible value of $k$ in a given problem. Tables \ref{tab:uniform_2} -- \ref{tab:uniform_mixed} display the grid convergence for several examples that were derived from known test solutions. Note that for each of these tables, changing the test solution only affects the source and boundary data, which means $\overline{Q}_N$ from \eqref{eq:N_domain_complete} and its QR-factorization remain the same  for all three test problems. Hence, for each of Tables \ref{tab:uniform_2} -- \ref{tab:uniform_mixed}, the QR-factorization is only computed once (during the first case) and is reused for the second and third cases, allowing those problems to be solved at a reduced cost.

    Table \ref{tab:uniform_2} shows the grid convergence of the case with two subdomains, which uses the example domain given in Figure~\ref{fig:simple_domains} with basic Robin boundary conditions defined by $\alpha = 1$ and $\beta = 1$ in \eqref{eq:helmholtz_robin}. By comparing the grid convergence of the first test case in Table \ref{tab:uniform_2} to that of Table \ref{tab:base_case}, we can see that including the domain decomposition does not affect the convergence rate of the method, and also has very little effect on the error itself. The second and third test solutions are both inhomogeneous, and clearly still converge with the designed rate of convergence.

    \begin{table}[ht]
        \centering
        \begin{tabular}{rcllcllcll}
        \toprule
            && \multicolumn{2}{c}{$e^{i\frac{k}{\sqrt{2}} (x+y)}$} && \multicolumn{2}{c}{$e^{\frac{-1}{1 - (x^2+y^2)}}$} && \multicolumn{2}{c}{$\sin^4(\pi x) \sin(\pi y)$}\\ \cmidrule{3-4} \cmidrule{6-7} \cmidrule{9-10}
            $n$ &&Error &Rate &&Error &Rate &&Error &Rate \\ \midrule
            64 &&8.52e-03 &- &&7.93e-04$^*$ &- &&1.37e-03 &- \\
            128 &&5.15e-04 &4.05 &&3.16e-06$^*$ &7.97 &&8.32e-05 &4.04 \\
            256 &&3.12e-05 &4.05 &&2.18e-07 &3.86 &&5.13e-06 &4.02 \\
            512 &&1.91e-06 &4.03 &&1.34e-08 &4.03 &&3.19e-07 &4.01 \\
            1024 &&1.20e-07 &4.00 &&8.31e-10 &4.01 &&2.00e-08 &4.00 \\
            2048 &&7.37e-09 &4.02 &&5.18e-11 &4.00 &&1.26e-09 &3.99 \\
        \bottomrule
        \end{tabular}
        \caption{Grid convergence for three test solutions over a two-subdomain $\Omega$ where $\Omega_1$ is centered at the origin and the interface with $\Omega_2$ is at $x=1$. A uniform wavenumber $k=13$ is used in both subdomains, $M^* = 40$, and the Robin boundary condition is uniformly defined by $\alpha = \beta = 1$ in \eqref{eq:helmholtz_robin}. Errors marked with $^*$ are computed with $M^*=20$.}
        \label{tab:uniform_2}
    \end{table}

    In Table \ref{tab:uniform_duct}, we use the exact same test solutions, wavenumber, and boundary conditions as in Table \ref{tab:uniform_2}, but on a larger scale with $N = 24$ subdomains extending in the positive $x-$direction from $\Omega_1$. On this larger scale, we can directly observe the $\mathcal{O} (N^3 n)$ complexity of the QR-factorization with respect to $n$, which was explained in Section~\ref{sec:complexity} (see Table \ref{tab:uniform_duct_times} for the complexity with respect to $N$). In the first two columns of Table \ref{tab:uniform_duct}, we see that as the grid dimension $n$ is doubled, the time for the QR-factorization is approximately doubled as well.

    \begin{table}[ht!]
        \centering
        \begin{tabular}{r c rr c ll c ll c ll}
        \toprule
            &&&& &\multicolumn{2}{c}{$e^{i\frac{k}{\sqrt{2}} (x+y)}$} &&\multicolumn{2}{c}{$e^{\frac{-1}{1 - (x^2+y^2)}}$} &&\multicolumn{2}{c}{$\sin^4(\pi x) \sin(\pi y)$}\\ \cmidrule{6-7} \cmidrule{9-10} \cmidrule{12-13}
            $n$ &&QR Time &Ratio &&Error &Rate &&Error &Rate &&Error &Rate \\ \midrule
            64 &&3.23 &- &&3.89e-02 &- &&8.34e-04$^*$ &- &&1.07e-03 &- \\
            128 &&6.58 &2.00 &&8.08e-04 &5.59 &&4.29e-06$^*$ &7.60 &&6.38e-05 &4.07 \\
            256 &&13.50 &2.05 &&4.85e-05 &4.06 &&6.10e-07 &2.81 &&3.95e-06 &4.01 \\
            512 &&28.91 &2.14 &&3.01e-06 &4.01 &&3.80e-08 &4.00 &&2.46e-07 &4.01 \\
            1024 &&60.13 &2.08 &&1.89e-07 &4.00 &&2.38e-09 &4.00 &&1.54e-08 &4.00 \\
            2048 &&108.77 &1.81 &&1.39e-08 &3.77 &&1.51e-10 &3.98 &&1.02e-09 &3.91 \\
        \bottomrule
        \end{tabular}
        \caption{Grid convergence and QR-factorization timing for three test solutions where $\Omega$ is a long duct of $N = 24$ subdomains. $\Omega_1$ is centered at the origin and each subsequent $\Omega_i$ is attached horizontally in the positive $x-$direction. A uniform wavenumber $k=13$ is used in all subdomains, $M^* = 40$, and the Robin boundary condition is uniformly defined by $\alpha = \beta = 1$ in \eqref{eq:helmholtz_robin}. The ratios of QR times demonstrates linear complexity with respect to the grid dimension $n$. Errors marked with $^*$ were computed with $M^*=20$. Note that the QR-factorization does not need repeated for the $^*$ cases, as their factorizations can be extracted directly from the existing factorization in each case.}
        \label{tab:uniform_duct}
    \end{table}

    \begin{table}[H]
        \centering
        \begin{tabular}{r c ll c ll c ll}
        \toprule
            &&\multicolumn{2}{c}{$e^{i\frac{k}{\sqrt{2}} (x+y)}$} &&\multicolumn{2}{c}{$e^{\frac{-1}{1 - (x^2+y^2)}}$} &&\multicolumn{2}{c}{$\sin^4(\pi x) \sin(\pi y)$}\\ \cmidrule{3-4} \cmidrule{6-7} \cmidrule{9-10}
            $n$ &&Error &Rate &&Error &Rate &&Error &Rate \\ \midrule
            64 &&9.08e-03 &- &&8.63e-04$^*$ &- &&1.86e-03 &- \\
            128 &&5.76e-04 &3.98 &&4.76e-06$^*$ &7.50 &&1.12e-04 &4.06 \\
            256 &&3.63e-05 &3.99 &&1.64e-07 &4.86 &&6.92e-06 &4.01 \\
            512 &&2.27e-06 &4.00 &&1.02e-08 &4.01 &&4.30e-07 &4.01 \\
            1024 &&1.42e-07 &4.00 &&6.39e-10 &4.00 &&2.70e-08 &4.00 \\
            2048 &&8.83e-09 &4.01 &&4.01e-11 &3.99 &&1.69e-09 &3.99 \\
        \bottomrule
        \end{tabular}
        \caption{Grid convergence for three test solutions where $\Omega$ is a large square comprising of 3 subdomains in each direction, and the bottom left subdomain is centered at the origin. A uniform wavenumber $k=13$ is used in all subdomains, $M^* = 40$, and the Robin boundary condition is uniformly defined by $\alpha = \beta = 1$ in \eqref{eq:helmholtz_robin}. Errors marked with $^*$ are computed with $M^* = 20$.}
        \label{tab:uniform_square}
    \end{table}

    The case of a large square domain is reported in Table \ref{tab:uniform_square}, with the same three test solutions as previous tables. In this case, there is a subdomain that is entirely interior and therefore has no boundary condition given, as well as numerous cross-points where more than two subdomains meet at a single point. As can be seen in Table \ref{tab:uniform_square}, the errors and convergence rate are unaffected by the presence of an interior subdomain and cross-points in this simple case, while more extreme cases are presented in Tables \ref{tab:self_square} and \ref{tab:piecewise_unique_k} in Section~\ref{sec:piecewise_k}.

    \begin{table}[ht]
        \centering
        \begin{tabular}{r c ll c ll c ll}
        \toprule
            &&\multicolumn{2}{c}{$e^{i\frac{k}{\sqrt{2}} (x+y)}$} &&\multicolumn{2}{c}{$e^{\frac{-1}{1 - (x^2+y^2)}}$} &&\multicolumn{2}{c}{$\sin^4(\pi x) \sin(\pi y)$}\\ \cmidrule{3-4} \cmidrule{6-7} \cmidrule{9-10}
            $n$ &&Error &Rate &&Error &Rate &&Error &Rate \\ \midrule
            64 &&1.73e-03 &- &&2.08e-04$^*$ &- &&9.61e-04 &- \\
            128 &&1.05e-04 &4.04 &&3.49e-06$^*$ &5.90 &&5.86e-05 &4.04 \\
            256 &&6.55e-06 &4.01 &&1.56e-07 &4.49 &&3.62e-06 &4.02 \\
            512 &&4.10e-07 &4.00 &&9.69e-09 &4.01 &&2.25e-07 &4.01 \\
            1024 &&2.56e-08 &4.00 &&6.06e-10 &4.00 &&1.41e-08 &3.99 \\
            2048 &&1.60e-09 &4.00 &&3.79e-11 &4.00 &&8.77e-10 &4.01 \\
        \bottomrule
        \end{tabular}
        \caption{Grid convergence for three test solutions over two subdomains with mixed boundary conditions as indicated in Figure~\ref{fig:uniform_mixed}. A uniform wavenumber $k=13$ is used in both subdomains, and $M^* = 40$. Errors marked with $^*$ are computed with $M^*=20$.}
        \label{tab:uniform_mixed}
    \end{table}
    
\begin{figure}[ht]
        \centering
        \includegraphics[width=.5\textwidth]{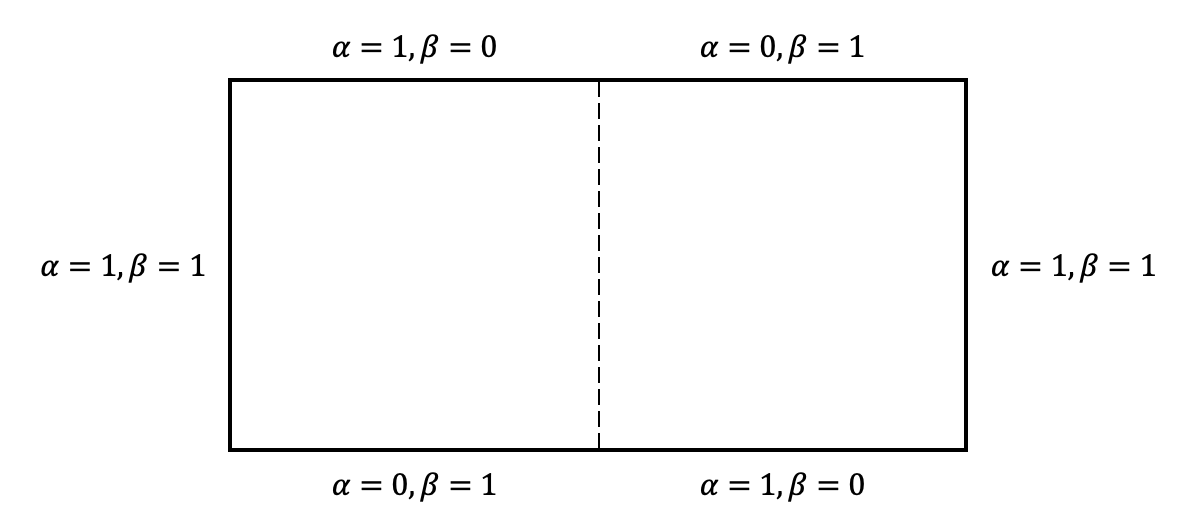}
        \caption{The mixed boundary condition for Table \ref{tab:uniform_mixed}. The coefficients $\alpha$ and $\beta$ are defined separately for each exterior edge of the subdomains.}
        \label{fig:uniform_mixed}
    \end{figure}

    Table \ref{tab:uniform_mixed} displays the grid convergence for a problem with mixed boundary conditions (see also Figure~\ref{fig:uniform_mixed}), showing that the method is robust enough to handle such boundary conditions without affecting its convergence rate. The domain, wavenumber, and test solutions in Table \ref{tab:uniform_mixed} are the same as in Table \ref{tab:uniform_2}, saving the cost of two applications of $\Gh$ (per subdomain) because $F$ in the right-hand side of \eqref{eq:N_block_system} does not need to be computed. As the type of boundary condition changed (i.e., $\alpha$ or $\beta$ in \eqref{eq:helmholtz_robin}~changed) we need to recompute the QR-factorization for the first test solution, reusing it for the second and third test solutions.

    In Table \ref{tab:uniform_duct_times}, we can see how the timing for the QR-factorization grows with respect to the number of subdomains, $N$. Recall from Section~\ref{sec:complexity} that the complexity of the QR-factorization should be $\mathcal{O} \left(N^3 n \right)$, so as $N$ is doubled we would expect to see the execution time of the QR-factorization increase by a factor of 8. However, as can be seen in Table \ref{tab:uniform_duct_times}, the execution time of the QR-factorization actually scales slower than its theoretical complexity would suggest, at least for the sizes of problems we were able to test.

    \begin{table}[ht]
        \centering
        \begin{tabular}{r c rr}
        \toprule
            $N$ &&QR time &Ratio \\ \midrule
            1 &&0.064 &- \\
            2 &&0.29 &4.61 \\
            4 &&1.66 &5.65 \\
            8 &&11.72 &7.04 \\
            16 &&61.89 &5.28 \\
            32 &&362.42 &5.86 \\
        \bottomrule
        \end{tabular}
        \caption{Timings for the QR-factorization (in seconds) on the $2048 \times 2048$ grid. $\Omega$ is a long duct of $N$ subdomains, and the test solution is $u = e^{i\frac{k}{\sqrt{2}} (x+y)}$ with $M^* = 40$ and $k = 13$ in all subdomains. The Robin boundary condition is defined by $\alpha = \beta = 1$ in \eqref{eq:helmholtz_robin}.}
        \label{tab:uniform_duct_times}
    \end{table}

    \subsection{Piecewise-Constant Wavenumber} \label{sec:piecewise_k}
    We now focus on cases where $k$ is piecewise-constant over $\Omega$ (with constant value $k_i$ over any $\Omega_i$). New $Q^{(i,*,*)}$ matrices are needed for any new values of $k_i$, but recall that we only need to compute $Q^{(i,*,*)}$ once for each unique $k_i$. For the results in this section, it is assumed that any necessary $Q^{(i,*,*)}$ matrices were appropriately computed ahead of time.

\begin{figure}[ht]
        \centering
        \begin{subfigure}{.31\textwidth}
            \centering
            \includegraphics[width=\textwidth]{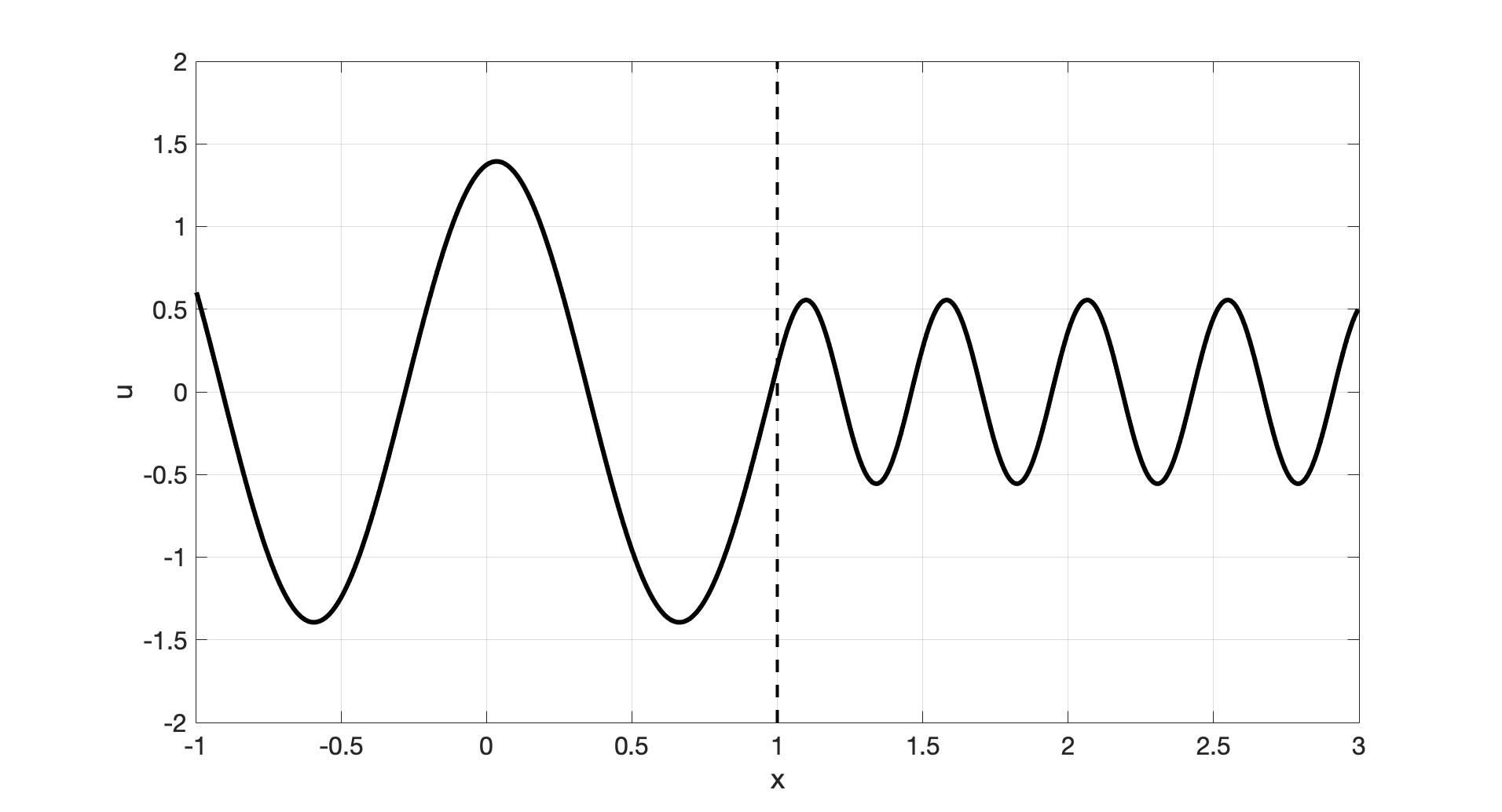}
            \caption{$k_1=5$, $k_2=13$}
            \label{fig:k13_jump}
        \end{subfigure}
        \begin{subfigure}{.31\textwidth}
            \centering
            \includegraphics[width=\textwidth]{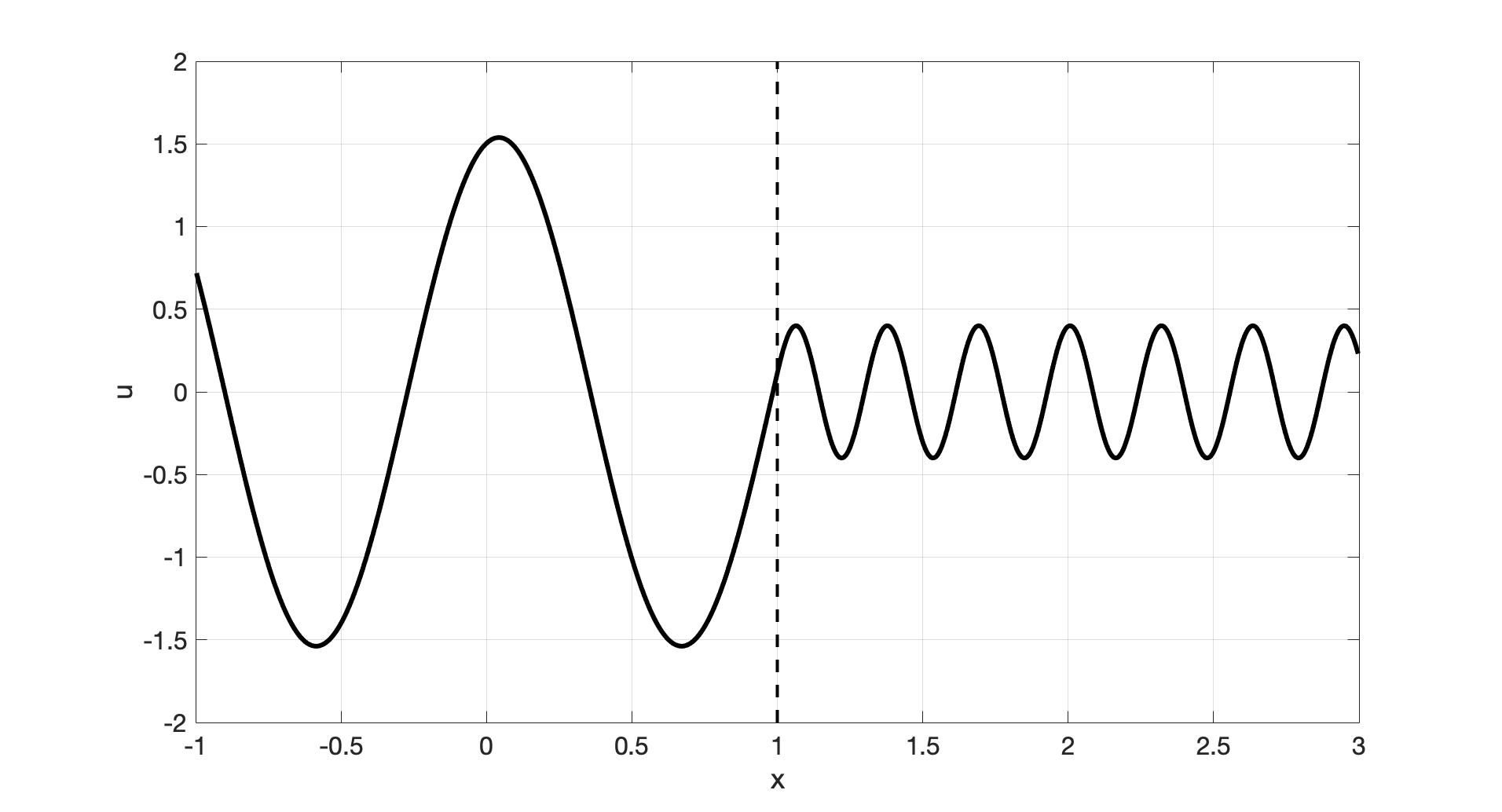}
            \caption{$k_1=5$, $k_2=20$}
            \label{fig:k20_jump}
        \end{subfigure}
        \begin{subfigure}{.31\textwidth}
            \centering
            \includegraphics[width=\textwidth]{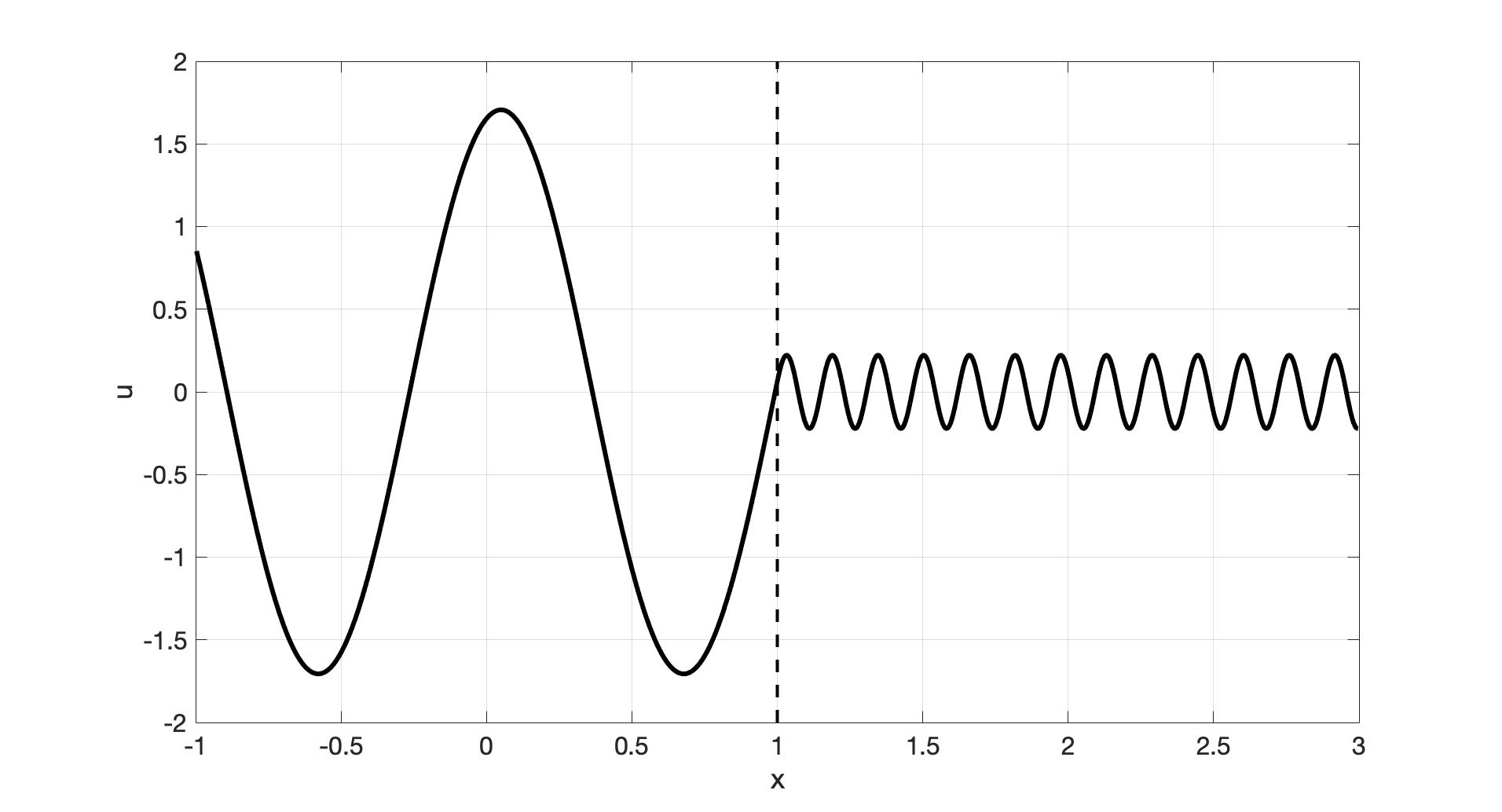}
            \caption{$k_1=5$, $k_2=40$}
            \label{fig:k40_jump}
        \end{subfigure}
        \caption{Real part of the test solutions from Table \ref{tab:piecewise_2}. These solutions have no dependence on $y$, so cross-sections in the $x-$direction are plotted. In each plot, it is clear that the frequency of the plane-wave changes at the interface.}
        \label{fig:kjump_plots}
    \end{figure}

\begin{figure}[ht]
        \centering
        \includegraphics[width=0.8\textwidth]{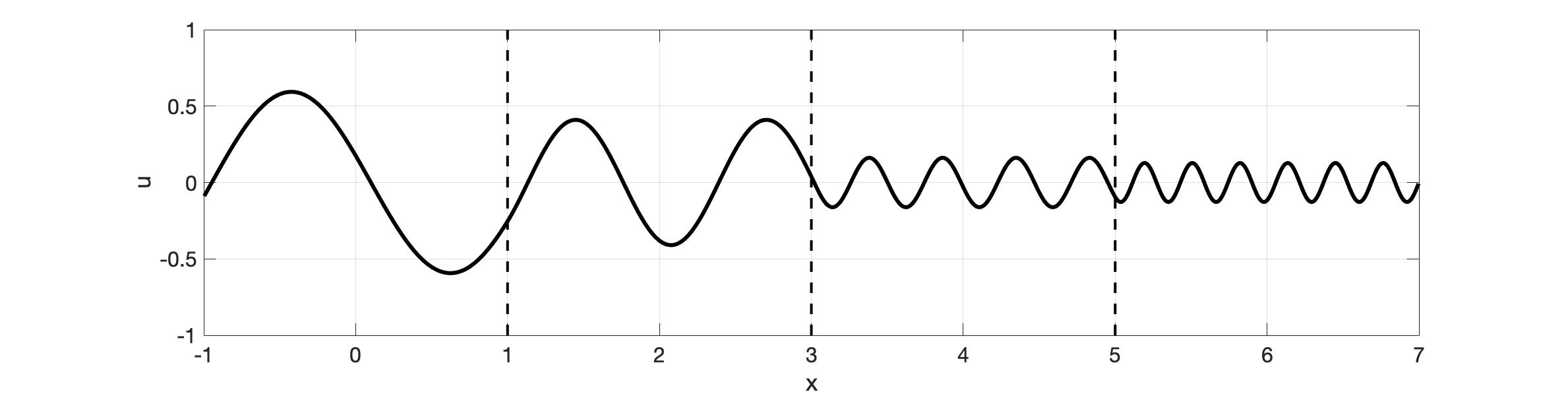}
        \caption{Real part of the test solution from Table \ref{tab:piecewise_4}. The solution has no dependence on $y$, so a cross-section is plotted. Each $\Omega_i$ has a distinct wavenumber $k_i$, with $k_1=3$, $k_2=5$, $k_3=13$, and $k_4=20$. Moving in the positive $x-$direction, the wavenumber and frequency increase, while the amplitude decreases.}
        \label{fig:piecewise_4}
    \end{figure}

    Tables \ref{tab:piecewise_2} and \ref{tab:piecewise_4} show the grid convergence for the two-subdomain and four-subdomain cases, respectively. In both cases, the test solution is obtained by considering an incident wave, $e^{ik_1x}$ in $\Omega_1$, and deriving the corresponding reflected and transmitted waves by enforcing continuity of the function and its normal derivative at each interface. This derivation can be found in Appendix \ref{app:piecewise}. Table \ref{tab:piecewise_2} shows the results of taking $k_1 \!=\! 5$ and allowing jumps to $k_2$, varying between $13$, $20$, and $40$, with cross-sections of each of these solutions plotted in Figure~\ref{fig:kjump_plots}. The method maintains its fourth-order rate of convergence, even on the largest jump from $5$ to $40$. It is worth pointing out here that as $k_2$ is increased, we need to increase $M^*$ because more oscillatory solutions will require more basis functions to maintain high-order accuracy. In cases where $M^* \geq 50$, this causes a loss of accuracy on the coarsest grid, so $M^*$ is reduced only for the $n=64$ grid in the relevant test cases. The four-subdomain case presented in Table \ref{tab:piecewise_4} was derived similar to the two-subdomain case, but only for one test solution. The values of $k$ on each subdomain in this example are $k_1=3$, $k_2=5$, $k_3=13$, and $k_4=20$, and the solution is plotted in Figure~\ref{fig:piecewise_4}. Even with four unique wavenumbers, it can be seen in Table \ref{tab:piecewise_4} that the method still has fourth-order convergence.

    \begin{table}[ht]
        \centering
        \begin{tabular}{r c ll c ll c ll}
        \toprule
            &&$k_1=5$ &$k_2=13$ &&$k_1=5$ &$k_2=20$ &&$k_1=5$ &$k_2=40$ \\ 
            &&\multicolumn{2}{c}{$M^*=40$} &&\multicolumn{2}{c}{$M^*=50$} &&\multicolumn{2}{c}{$M^*=60$}\\ \cmidrule{3-4} \cmidrule{6-7} \cmidrule{9-10}
            $n$ &&Error &Rate &&Error &Rate &&Error &Rate \\ \midrule
            64 &&3.23e-02 &- &&8.81e-02$^*$ &- &&1.27e+01$^*$ &- \\
            128 &&1.98e-03 &4.03 &&5.13e-03 &4.10 &&1.08e-01 &6.88\\
            256 &&1.21e-04 &4.04 &&3.13e-04 &4.03 &&6.51e-03 &4.06\\
            512 &&7.59e-06 &3.99 &&1.96e-05 &4.00 &&4.03e-04 &4.01\\
            1024 &&4.69e-07 &4.02 &&1.22e-06 &4.01 &&2.52e-05 &4.00\\
            2048 &&3.09e-08 &3.92 &&7.78e-08 &3.97 &&1.57e-06 &4.01\\
        \bottomrule
        \end{tabular}
        \caption{Grid convergence for the two-subdomain test case with the incident wave $u=e^{ik_1x}$ and Dirichlet boundary conditions, as plotted in Figure~\ref{fig:kjump_plots}. The jump in wavenumber goes from $k_1 = 5$ to the indicated value of $k_2$, and $M^*$ is chosen separately for each case to ensure accuracy beyond that obtained on the finest grid. Errors marked with a $^*$ were computed with $M^* = 30$.}
        \label{tab:piecewise_2}
    \end{table}

    \begin{table}[ht]
        \centering
        \begin{tabular}{r c l l}
        \toprule
            $n$ &&Error &Rate \\ \midrule
            64 &&2.02e-02$^*$ &- \\
            128 &&1.21e-03 &4.08 \\
            256 &&7.38e-05 &4.03 \\
            512 &&4.61e-06 &4.00 \\
            1024 &&2.87e-07 &4.01 \\
            2048 &&2.05e-08 &3.80 \\
        \bottomrule
        \end{tabular}
        \caption{Grid convergence for $u = e^{ikx}$ over four subdomains with Dirichlet boundary conditions, and wavenumbers $k_1=3$, $k_2=5$, $k_3=13$, and $k_4=20$. $\Omega_1$ is centered at the origin, and each subdomain extends in the positive $x-$direction, with $M^* = 50$. The error marked with a $^*$ was computed with $M^* = 30$. }
        \label{tab:piecewise_4}
    \end{table}

    Additionally, we point out the increase in error as $k_2$ increases in Table \ref{tab:piecewise_2}. We attribute this increase to the pollution effect \cite{acc}, because it appears consistent with our additional observations of the pollution effect for problems with uniform wavenumbers (no jumps) that are comparable to $k_2$. We therefore conclude that the method is not inherently sensitive to discontinuities in the wavenumber.

\begin{figure}[ht]
        \centering
        \includegraphics[width=0.8\textwidth]{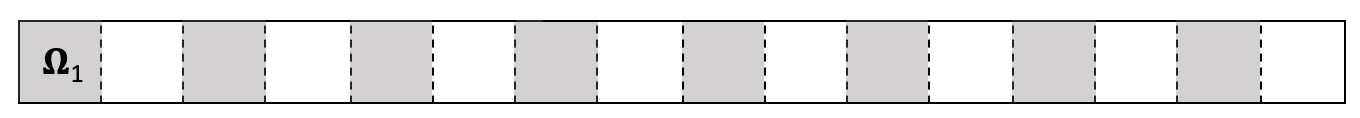}
        \caption{The decomposition used to compute Table \ref{tab:self_alternatingduct} where the wavenumbers are alternating in each subdomain. $\Omega_1$ is indicated in the left-most subdomain with each subsequent subdomain being attached in the positive $x-$direction. Wavenumbers are assigned as $k=5$ for gray subdomains and $k=40$ in white subdomains.}
        \label{fig:alternating_duct}
    \end{figure}

    As we allow the test cases to become more complex in geometry and wavenumber distribution, it becomes more difficult to obtain analytic test solutions. Instead, we specify the boundary and source data directly, and calculate errors on shared nodes between subsequent resolutions of the grid. For simplicity, the source function is taken to be a ``bump'' function:
    \begin{equation} \label{eq:source_bump}
        f(x,y) = \begin{cases} \exp{\left( \frac{-1}{\frac{1}{4} - (x^2 + y^2)} \right)} &x^2 + y^2 < \frac{1}{2} \\ 0 &\text{otherwise} \end{cases}
    \end{equation}
    In Table \ref{tab:self_alternatingduct}, we present an example of a long duct of $N = 16$ subdomains with a change in wavenumber at every interface, alternating between $k = 5$ and $k = 40$ (depicted in Figure~\ref{fig:alternating_duct}). The three cases in Table \ref{tab:self_alternatingduct} represent homogeneous Dirichlet, Neumann, and Robin ($\alpha=\beta=1$) boundary conditions, respectively, and show that the method maintains its design rate of convergence for all three types of boundary conditions.

    \begin{table}[ht]
        \centering
        \begin{tabular}{r c ll c ll c ll}
        \toprule
            &&\multicolumn{2}{c}{Dirichlet} &&\multicolumn{2}{c}{Neumann} &&\multicolumn{2}{c}{Robin}\\ \cmidrule{3-4} \cmidrule{6-7} \cmidrule{9-10}
            $n$ &&\selfnorm &Rate &&\selfnorm &Rate &&\selfnorm &Rate\\ \midrule
            256 &&1.51e-03 &- &&5.79e-04 &- &&3.50e-04 &- \\
            512 &&1.08e-04 &3.81 &&7.54e-06 &6.26 &&3.57e-06 &6.62 \\
            1024 &&6.76e-06 &3.99 &&4.62e-07 &4.03 &&2.19e-07 &4.03 \\
            2048 &&4.22e-07 &4.00 &&2.90e-08 &3.99 &&1.38e-08 &3.99 \\
        \bottomrule
        \end{tabular}
        \caption{Grid convergence on the duct of $N=16$ subdomains depicted in Figure~\ref{fig:alternating_duct}, various homogeneous boundary conditions, and the source function from \eqref{eq:source_bump}. The wavenumbers alternate between $k = 5$ and $k = 40$, and $M^* = 60$.}
        \label{tab:self_alternatingduct}
    \end{table}

\begin{figure}[ht!]
        \centering
        \includegraphics[width=0.2\textwidth]{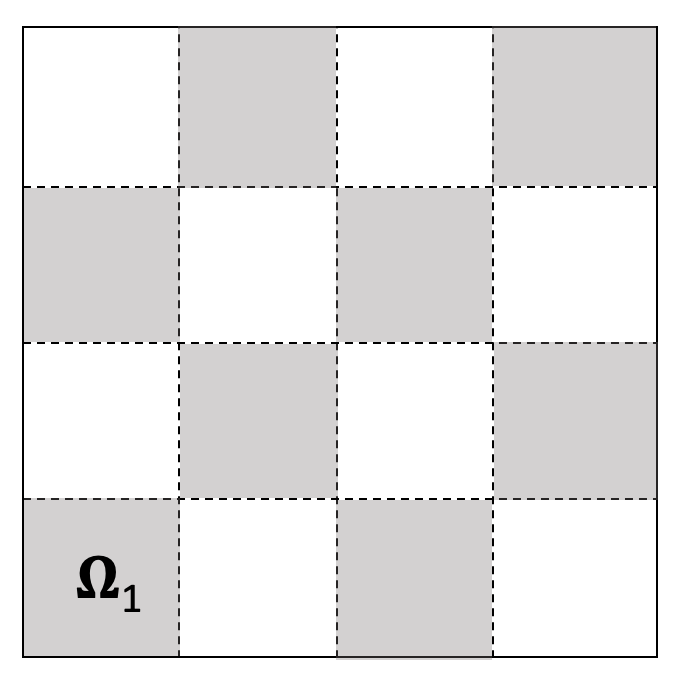}
        \caption{The $4 \times 4$ decomposition used in the first case of Table \ref{tab:self_square} where the wavenumbers are assigned in a checkerboard pattern, with $k=5$ in the gray subdomains and $k=40$ in the white subdomains.}
        \label{fig:checkerboard_square}
    \end{figure}

    In Table \ref{tab:self_square}, the domain is a large square decomposed into $N = N_d^2$ (cf. Section~\ref{sec:complexity}) subdomains, where the piecewise constant values of the wavenumber are defined in a checkerboard pattern with $k=5$ and $k=40$ as in Figure~\ref{fig:checkerboard_square}. These examples combine the qualitative aspects of Tables \ref{tab:uniform_square} and \ref{tab:self_alternatingduct}, containing cross-points as well as a changing wavenumber at every interface (now in both the $x-$ and $y-$directions). We emphasize  that no special considerations were given to internal or boundary cross-points, yet the method's convergence does not suffer from their presence.

    \begin{table}[ht]
        \centering
        \begin{tabular}{r c cc c cc c cc}
        \toprule
            &&\multicolumn{2}{c}{$N=16$} &&\multicolumn{2}{c}{$N=25$} &&\multicolumn{2}{c}{$N=36$} \\ \cmidrule{3-4} \cmidrule{6-7} \cmidrule{9-10}
            $n$ &&\selfnorm &Rate &&\selfnorm &Rate &&\selfnorm &Rate \\ \midrule
            256 &&8.54e-04 &- &&6.48e-03 &- &&4.31e-04 &-  \\
            512 &&6.90e-05 &3.63 &&7.45e-04 &3.12 &&1.76e-05 &4.61  \\
            1024 &&4.36e-06 &3.98 &&4.90e-05 &3.93 &&1.09e-06 &4.01  \\
            2048 &&2.73e-07 &4.00 &&3.06e-06 &4.00 &&6.79e-08 &4.01  \\
        \bottomrule
        \end{tabular}
        \caption{Grid convergence on a large square decomposed into $N=N_d^2$ subdomains and a checkerboard pattern for its wavenumber as depicted in Figure~\ref{fig:checkerboard_square}. There are homogeneous Dirichlet boundary conditions and the source function is \eqref{eq:source_bump}. The checkerboard wavenumbers are $k=5$ (gray) and $k=40$ (white), with $M^* = 60$.}
        \label{tab:self_square}
    \end{table}

\begin{figure}
        \centering
        \includegraphics[scale=0.4]{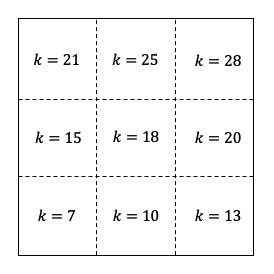}
        \caption{Piecewise constant values of the wavenumber $k$ for the example computed in Table \ref{tab:piecewise_unique_k}.}
        \label{fig:k_profile}
    \end{figure}

    Table \ref{tab:piecewise_unique_k} shows the final example, which returns to the case of a square domain decomposed into $3\times 3$ subdomains, but with wavenumbers assigned as in Figure~\ref{fig:k_profile}. In contrast to the configurations of Figures \ref{fig:alternating_duct} and \ref{fig:checkerboard_square}, this example uses a different wavenumber in each subdomain (similar to Figure~\ref{fig:piecewise_4}). This is the costliest test case for the method, as $Q^{(i,*,*)}$ is different for every $\Omega_i$. However, we can still use those $Q^{(i,*,*)}$ matrices computed for $k=5,13$, and $20$ that were used for earlier examples. In Table \ref{tab:piecewise_unique_k}, fourth-order convergence is  clearly observed.
    \begin{table}[ht]
        \centering
        \begin{tabular}{r c cc}
        \toprule
            $n$ &&\selfnorm &Rate \\ \midrule
            256 &&9.04e-05 &- \\
            512 &&3.14e-06 &4.85 \\
            1024 &&1.89e-07 &4.05 \\
            2048 &&1.20e-08 &3.98 \\
        \bottomrule
        \end{tabular}
        \caption{Grid convergence on a square domain decomposed into $N=9$ subdomains with homogeneous Dirichlet boundary conditions and the source function from \eqref{eq:source_bump}. Each subdomain has a unique wavenumber as depicted in Figure~\ref{fig:k_profile}, and $M^* = 60$.}
        \label{tab:piecewise_unique_k}
    \end{table}

    \subsection{Further Studies} \label{sec:further}
\begin{figure}
        \centering
        \includegraphics[scale=.3]{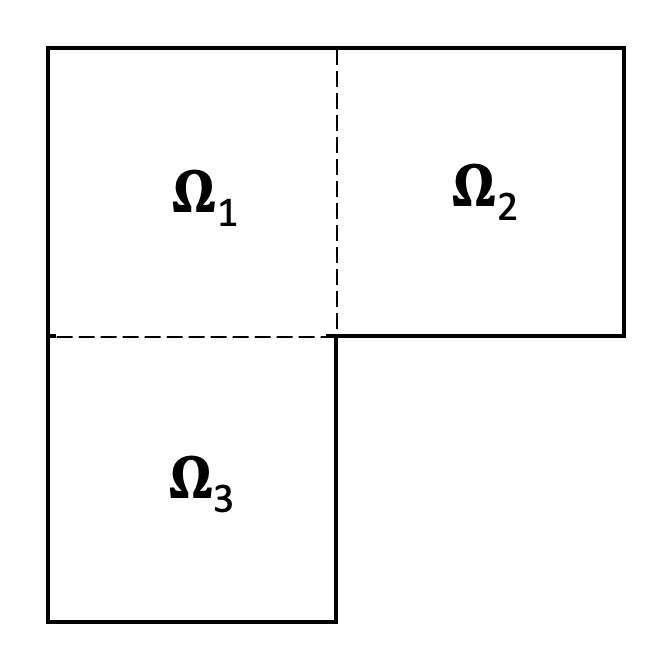}
        \caption{The ``Block L'' domain decomposition used in Table \ref{tab:blockL_convergence}.  $\Omega_1$ is centered over the origin, with interfaces to $\Omega_2$ and $\Omega_3$ at $x=1$ and $y=-1$, respectively. This configuration creates a reentrant corner at the point $(x,y) = (1,-1)$.}
        \label{fig:blockL_domain}
    \end{figure}

    Thus far, we have focused on domains that are either squares or ducts, making it straightforward to generate solutions that have no singularities and can ensure the design rate of convergence for the method. Moving beyond simple square and duct decompositions will often introduce reentrant corners, an example of which is shown in the third plot of Figure~\ref{fig:N_domains_valid}. Reentrant corners may cause a solution to develop a singularity on the boundary, which in turn can cause the method to lose its design rate of convergence. In order to observe this phenomenon, we introduce the ``Block L'' domain in Figure~\ref{fig:blockL_domain} which contains one reentrant corner, and demonstrate the convergence obtained by our method under two different types of solutions.

    Table \ref{tab:blockL_convergence} shows the grid convergence over the Block L domain. The first case is derived from the exact solution $u = e^{i\frac{k}{\sqrt{2}}(x+y)}$, and the second case uses a homogeneous Dirichlet boundary condition with the source function \eqref{eq:source_bump}. In the first case, the test solution contains no singularities, and as expected Table \ref{tab:blockL_convergence} reflects the design rate of convergence. This indicates that the reentrant corner itself is not an issue for our method. However, the solution in the second case develops a singularity, and the breakdown in the convergence rate in Table \ref{tab:blockL_convergence} is indicative of that. For Case 1, we are able to compute the error directly with the known test solution, but for Case 2 we compute the error on successively refined grids as described earlier.

    The breakdown of convergence in the second case is a result of the solution's own singularity, rather than a shortcoming of the method. The resolution of singularities at reentrant corners with the MDP has been explored in \cite{Magura2017}. The method can also handle general shaped subdomains but these fall outside the scope of this paper.

    \begin{table}[ht]
        \centering
        \begin{tabular}{r c cc c cc}
        \toprule
            &&\multicolumn{2}{c}{Case 1} &&\multicolumn{2}{c}{Case 2} \\ \cmidrule{3-4} \cmidrule{6-7}
            $n$ &&Error &Rate &&\selfnorm &Rate \\ \midrule
            64 &&9.58e-04 &- &&- &- \\
            128 &&6.03e-05 &3.99 &&4.08e-06 &- \\
            256 &&3.76e-06 &4.00 &&1.98e-06 &1.04 \\
            512 &&2.35e-07 &4.00 &&1.09e-06 &0.86 \\
            1024 &&1.47e-08 &4.00 &&1.02e-06 &0.10 \\
            2048 &&9.80e-10 &3.91 &&4.70e-07 &1.12 \\
        \bottomrule
        \end{tabular}
        \caption{Grid convergence for the Block L case from Figure~\ref{fig:blockL_domain} with Dirichlet boundary conditions and uniform wavenumber $k=13$, with $M^*=40$. In Case 1, the boundary and source data are derived from the plane-wave $u = e^{i\frac{k}{\sqrt{2}}(x+y)}$, which results in fourth order convergence. In Case 2, the boundary data is zero, and the source function is given in \eqref{eq:source_bump}. Case 2 develops a singularity at the reentrant corner and breaks down the convergence.}
        \label{tab:blockL_convergence}
    \end{table}

\section{Conclusions} \label{sec:conclusions}
In this paper, we have adapted the Method of Difference Potentials to solve a non-overlapping Domain Decomposition formulation for the Helmholtz equation. After solving for the boundary information along all interfaces, the direct solves for all subproblems can be distributed and performed concurrently. Further, the formulation is convenient for handling piecewise-constant wavenumbers, as well as mixed boundary conditions. Numerical results corroborate the fourth-order convergence rate of the method in numerous situations, most notably for decompositions with cross-points and for transmission problems with a large jump in the wavenumber. Our formulation also allowed us to demonstrate different behaviors of the method, such as its performance on solutions with singularities and the method's complexity with respect to the number of subdomains or grid dimension.

Once the boundary/interface data have been obtained for the subdomains, only one direct solve is required per subdomain, and this set of direct solves can be parallelized on a number of processors up to the number of subdomains. However, this comes at the cost of requiring the QR-factorization of a large system. Even though the QR-factorization scales slower than its theoretical complexity would suggest, for cases with more than $N=8$ subdomains the QR-factorization already outweighs the final PDE solves in cost. However, once the QR-factorization has been computed, new problems with variations in the source and boundary data can be solved at a reduced cost by simply reusing the computed QR-factorization. The framework that is laid out in this paper can be adapted and extended in numerous ways to broaden the applicability of the method. Briefly mentioned in Section~\ref{sec:transmission_condition}, transmission conditions of the form \eqref{eq:transmission_general} can be implemented trivially along every interface to account for more complex properties of the solution there, such as jumps in the normal derivative or  solution itself. The method can be generalized to account for a smoothly varying wavenumber $k = k(x,y)$ in each subdomain, although it may reduce the efficiency  as the FFT-based solver would no longer be applicable. Base subdomains of a different shape could  be included, such as those with piecewise-curvilinear boundaries, in order to account for more complex geometries.

%% The Appendices part is started with the command \appendix;
%% appendix sections are then done as normal sections
\appendix
\section{Derivation of a Function with Piecewise Constant Wavenumber} \label{app:piecewise}
\begin{figure}[ht]
    \centering
    \includegraphics[width=0.4\textwidth]{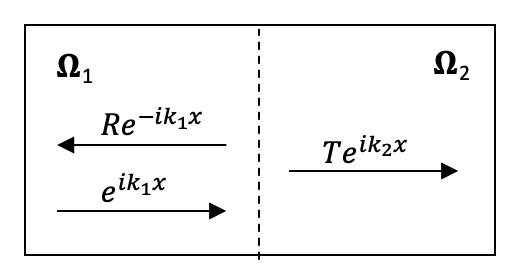}
    \caption{The setup for deriving the reflected and transmitted parts of a one-dimensional incident wave hitting an interface in a domain composed of two subdomains.}
    \label{fig:two_derivation}
\end{figure}

Consider a domain $\Omega \subset \R^2$ split into two subdomains, $\Omega_1$ and $\Omega_2$, as in Figure~\ref{fig:two_derivation}. Let each subdomain have its own corresponding wavenumber, $k_1$ or $k_2$, which we will assume is constant for simplicity in this derivation. We seek a function $u \in C^1(\Omega)$ that incorporates the reflected and transmitted parts of an incident wave that starts in $\Omega_1$ and propagates toward the interface, located at $x=0$ for simplicity.

Let $u_1 = e^{ik_1x}$ be the incident wave in $\Omega_1$. When the wavenumber changes at the interface $x=0$, $u_1$ is partially reflected back into $\Omega_1$ and partially transmitted through to $\Omega_2$. The reflected wave has some amplitude $R$ and travels in the opposite direction of $u_1$, giving us $u_2 = Re^{-ik_1x}$. The transmitted part, on the other hand, will have its own amplitude $T$, traveling in the same direction as $u_1$ and with the wavenumber $k_2$, giving $u_3 = Te^{ik_2x}$. This allows us to write the function $u$ as:
\begin{equation} \label{eq:cases_function}
    u(x,y) = \begin{cases} e^{ik_1x} + Re^{-ik_1x}, &x \leqslant 0 \\ Te^{ik_2x}, &x\geqslant 0 \end{cases}
\end{equation}
The condition to enforce continuity of the function at the interface is
\begin{align} \label{eq:condition_function_2}
    e^{ik_1x}|_{x=0} + Re^{-ik_1x}|_{x=0} &= Te^{ik_2x}|_{x=0} \nonumber \\
    \implies 1 + R &= T
\end{align}
and for continuity of the derivative we get
\begin{align} \label{eq:condition_derivative_2}
    \left( \partialdiff{}{x}e^{ik_1x} \right) |_{x=0} + \left( \partialdiff{}{x}Re^{-ik_1x} \right) |_{x=0} &= \left( \partialdiff{}{x}Te^{ik_2x} \right) |_{x=0} \nonumber \\
    \implies \left( ik_1 e^{ik_1x} \right) |_{x=0} - \left( ik_1 Re^{-ik_1x} \right) |_{x=0} &= \left( ik_2 Te^{ik_2x} \right) |_{x=0} \nonumber \\
    \implies ik_1 - ik_1 R &= ik_2 T \nonumber \\
    \implies k_1 (1 - R) &= k_2 T \nonumber \\
    \implies 1 - R &= \frac{k_2}{k_1} T
\end{align}
By combining \eqref{eq:condition_function_2} and \eqref{eq:condition_derivative_2}, we can solve for $R$ and $T$ to get
\begin{equation} \label{eq:R_T_solved}
    R = \frac{k_2}{2k_1} - \frac{1}{2} \hspace{2cm} T = \frac{k_2}{2k_1} + \frac{1}{2}
\end{equation}
The values of $R$ and $T$ from \eqref{eq:R_T_solved} can be plugged into \eqref{eq:cases_function} to obtain the function $u(x,y)$, defined across $\Omega$.

\begin{figure}[ht]
    \centering
    \includegraphics[width=0.6\textwidth]{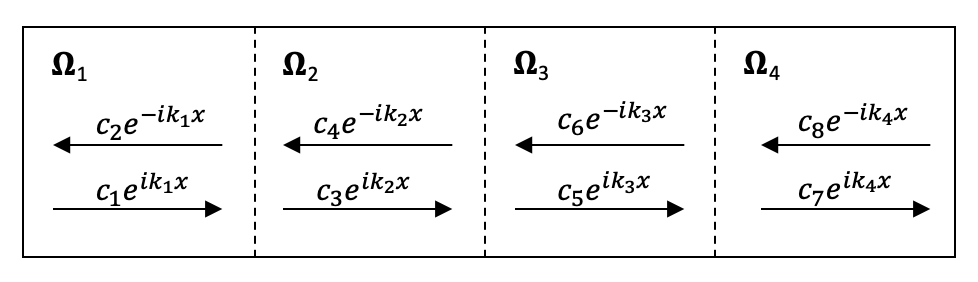}
    \caption{The setup for deriving the reflected and transmitted parts of a one-dimensional incident wave hitting an interface in a domain composed of four subdomains.}
    \label{fig:four_derivation}
\end{figure}

For a larger case with four subdomains, consider the scenario depicted in Figure~\ref{fig:four_derivation}, with interfaces at $x=x_1,\,x_2,$ and $x_3$. This scenario is a direct extension of the two-subdomain case, and we can obtain a linear system by similarly enforcing continuity of the function and its derivative at each interface. For example, at the interface between $\Omega_1$ and $\Omega_2$ (i.e. $x = x_1$) enforcing continuity of the function itself yields:
\begin{equation*}
    c_1e^{ik_1x_1} + c_2e^{-ik_1x_1} = c_3e^{ik_2x_1} + c_4e^{-ik_2x_1}
\end{equation*}
which can be rewritten as
\begin{equation*}
    c_1e^{ik_1x_1} + c_2e^{-ik_1x_1} - c_3e^{ik_2x_1} - c_4e^{-ik_2x_1} = 0
\end{equation*}
By including the corresponding conditions for both the function and its derivative at all three interfaces, we get the following system of equations:
\begin{align*}
    c_1e^{ik_1x_1} + c_2e^{-ik_1x_1} -c_3e^{ik_2x_1} - c_4e^{-ik_2x_1} &= 0 \\
    c_3e^{ik_2x_2} + c_4e^{-ik_2x_2} -c_5e^{ik_3x_2} - c_6e^{-ik_3x_2} &= 0 \\
    c_5e^{ik_3x_3} + c_6e^{-ik_3x_3} -c_7e^{ik_4x_3} - c_8e^{-ik_4x_3} &= 0 \\
    k_1c_1e^{ik_1x_1} - k_1c_2e^{-ik_1x_1} -k_2c_3e^{ik_2x_1} + k_2c_4e^{-ik_2x_1} &= 0 \\
    k_2c_3e^{ik_2x_2} - k_2c_4e^{-ik_2x_2} -k_3c_5e^{ik_3x_2} + k_3c_6e^{-ik_3x_2} &= 0 \\
    k_3c_5e^{ik_3x_3} - k_3c_6e^{-ik_3x_3} -k_4c_7e^{ik_4x_3} + k_4c_8e^{-ik_4x_3} &= 0
\end{align*}
Note that this only provides six equations for eight unknowns. As in the two-subdomain case, we can pick one of the waves to be the incident wave, and choose to set its amplitude to $1$ for simplicity, so we can directly choose $c_1 = 1$. Further, the boundary of $\Omega$ is not reflective, which means that $c_7e^{ik_4x}$ does not reflect upon reaching the right boundary, leaving $c_8 = 0$. This leaves six equations for six unknowns, which  allows this problem to be solved uniquely.

The same process can be applied to the $N-$subdomain case. If we let the incident wave to be given in $\Omega_1$, then its reflection back into $\Omega_1$ has an undetermined amplitude. For $\Omega_2$ through $\Omega_{N-1}$, there are two waves traveling in opposite directions for which the amplitudes are also undetermined. Finally, there is no reflected wave in $\Omega_N$, so there is only one amplitude to solve for, yielding a total of $1 + 2(N-2) + 1 = 2N-2$ unknowns. This scenario contains $N-1$ interfaces, and each interface has two conditions: continuity of the function and continuity of the derivative. These conditions yield $2(N-1) = 2N-2$ equations, allowing us to solve for the $2N-2$ unknowns.

%% References
%%
%% Following citation commands can be used in the body text:
%% Usage of \cite is as follows:
%%   \cite{key}          ==>>  [#]
%%   \cite[chap. 2]{key} ==>>  [#, chap. 2]
%%   \citet{key}         ==>>  Author [#]

%% References with bibTeX database:

% \bibliographystyle{model1-num-names}

%% New version of the num-names style
\bibliographystyle{elsarticle-num-names}
\bibliography{EvanNorth-elsevier.bib}

%% Authors are advised to submit their bibtex database files. They are
%% requested to list a bibtex style file in the manuscript if they do
%% not want to use model1-num-names.bst.

\end{document}